\documentclass[11pt]{article}
\usepackage{geometry}
\usepackage{authblk}
\usepackage{amsmath,amssymb,graphicx,tikz,multirow,hhline,color,array}
\usepackage[labelfont=bf]{caption}
\usepackage{mathtools}

\usepackage{natbib}
\bibpunct{[}{]}{,}{n}{}{}

\definecolor{mygray}{rgb}{.6, .6, .6}

\allowdisplaybreaks

\newcommand{\beq}{\begin{equation}}
\newcommand{\eeq}{\end{equation}}

\newcommand{\ov}{\overline}
\newcommand{\x}{\mathbf{x}}
\newcommand{\f}{\mathbf{f}}
\newcommand{\p}{\mathbf{p}}
\newcommand{\bb}{\mathbf{b}}
\newcommand{\Q}{\mathcal{Q}}
\newcommand{\bhat}{\mathbf{\hat b}}
\newcommand{\Qhat}{\mathbf{\hat{\mathcal{Q}}}}
\newcommand{\mom}{{x^ly^mz^n}}
\newcommand{\nmom}{\overline{\protect\raisebox{0pt}[10pt]{$\overline\mom$}}}

\title{Bounding averages rigorously using semidefinite~programming: mean moments of the Lorenz system}
\author{David Goluskin\thanks{Email: {\tt goluskin@uvic.ca}}}
\affil{Department of Mathematics and Statistics, University of Victoria\\
Victoria, BC, V8P 5C2, Canada and \\
Department of Mathematics, University of Michigan, Ann Arbor, MI 48109, USA}

\begin{document}

\maketitle

\begin{abstract}
We describe methods for proving bounds on infinite-time averages in differential dynamical systems. The methods rely on the construction of nonnegative polynomials with certain properties, similarly to the way nonlinear stability can be proved using Lyapunov functions. Nonnegativity is enforced by requiring the polynomials to be sums of squares, a condition which is then formulated as a semidefinite program (SDP) that can be solved computationally. Although such computations are subject to numerical error, we demonstrate two ways to obtain rigorous results: using interval arithmetic to control the error of an approximate SDP solution, and finding exact analytical solutions to relatively small SDPs. Previous formulations are extended to allow for bounds depending analytically on parametric variables. These methods are illustrated using the Lorenz equations, a system with three state variables $(x,y,z)$ and three parameters $(\beta,\sigma,r)$. Bounds are reported for infinite-time averages of all eighteen moments $x^ly^mz^n$ up to quartic degree that are symmetric under $(x,y)\mapsto(-x,-y)$. These bounds apply to all solutions regardless of stability, including chaotic trajectories, periodic orbits, and equilibrium points. The analytical approach yields two novel bounds that are sharp: the mean of $z^3$ can be no larger than its value of $(r-1)^3$ at the nonzero equilibria, and the mean of $xy^3$ must be nonnegative. The interval arithmetic approach is applied at the standard chaotic parameters to bound eleven average moments that all appear to be maximized on the shortest periodic orbit. Our best upper bound on each such average exceeds its value on the maximizing orbit by less than 1\%. Many bounds reported here are much tighter than would be possible without computer assistance.
\end{abstract}

\section{\label{sec: intro} Introduction}

In the study of dynamical systems, especially chaotic systems, time-averaged quantities are often of more interest than the details of any particular solution. Such quantities are typically estimated by numerically integrating the system and averaging along the resulting trajectory. Here we take the complementary approach of proving bounds directly on infinite-time averages. This approach offers several advantages. In a dynamical system with parametric variables, bounds can apply to entire parameter ranges and can depend analytically on the parameters, whereas numerical integration gives information only for the parameter values at which it is carried out. In a dynamical system with multiple local attractors, bounds can apply to all attractors, whereas a numerically integrated solution gives information only about the attractor in whose basin it begins. Furthermore, it is natural to derive bounds in a mathematically rigorous way.

The primary challenge in using bounds to estimate mean quantities is to construct bounds that are sufficiently tight, meaning they are close to the true values being bounded. Results that can be proved ``by hand'' are often quite conservative, but computer-assisted methods introduced by Chernyshenko \emph{et al.}\ \cite{Chernyshenko2014} and developed further by Fantuzzi \emph{et al.}\ \cite{Fantuzzi2016} have given tighter results in several examples. These methods are part of an approach used often in the study of dynamics: constructing functions with certain properties that imply the desired result. The most widely known example is that of Lyapunov functions, which imply nonlinear stability. Here we construct \emph{auxiliary functions} that have related but distinct properties and imply bounds on time averages. For all well-posed systems whose trajectories remain in a compact region, it has been proved that optimally chosen auxiliary functions imply arbitrarily sharp bounds \cite{Tobasco2017}. Much like Lyapunov's method, the main difficulty is to choose auxiliary functions optimally. In many cases this challenge can be made tractable by sum-of-squares (SOS) relaxation, wherein nonnegativity of a function is replaced with the stronger condition that the function is representable as a sum of squares of polynomials \cite{Parrilo2003}. Such relaxation is useful because an SOS condition can be equivalently formulated as a semidefinite program (SDP), a type of convex optimization problem solvable by various software packages \cite{Boyd2004}. The bounding methods developed in \cite{Chernyshenko2014, Fantuzzi2016} follow this same philosophy; SOS conditions are formulated that imply the desired bounds, and these SOS conditions are then reformulated as SDPs.

The approach of bounding time averages using SDPs was applied to the van der Pol oscillator in \cite{Fantuzzi2016}, and very tight bounds were obtained for averages over the limit cycle. This success is promising, but the phase space of the van der Pol system is simple, consisting of one limit cycle and one repelling equilibrium. In the present work we demonstrate that the SDP approach to bounding can succeed also for systems with much more complicated phase spaces. In particular we consider the Lorenz system \cite{Lorenz1963},
\begin{align}
\tfrac{d}{dt}x &= -\sigma x + \sigma y, &
\tfrac{d}{dt}y &= rx - y - xz, &
\tfrac{d}{dt}z &= xy - \beta z,
\label{eq: Lorenz}
\end{align}
and bound time averages of moments of the coordinates -- that is, monomials $\mom$ where $l,m,n$ are nonnegative integers. The bounds constructed here are global in the sense that they hold for all possible trajectories. An upper bound on a time average is sharp if it equals the maximum of that average over all trajectories; a lower bound is sharp if it equals the minimum over all trajectories. Methods for proving bounds holding only in particular basins of attraction (or in the presence of noise) are developed elsewhere \cite{Fantuzzi2016}. Some of our results apply to large sets of the parameters $(\beta,\sigma,r)$, and others are specific to the standard chaotic values $(8/3,10,28)$. We assume $\sigma\neq0$ throughout.

At the standard parameter values we have bounded time averages of eighteen moments -- all moments up to quartic degree that are invariant under the symmetry of the Lorenz system. All eighteen averages can be bounded below by zero, and these lower bounds are sharp since they are saturated by trajectories on or tending to the unstable equilibrium at the origin. For seven of the moments ($z$, $x^2$, $xy$, $z^2$, $xyz$, $z^3$, $xyz^2$) we construct sharp upper bounds, some of which have been proved before. In all seven cases the upper bounds are saturated by trajectories on or tending to either of the unstable nonzero equilibria. The same appears true for $x^2z$, but our upper bound is not quite sharp. Time averages of the remaining ten moments all seem to be maximized by trajectories on or tending to the shortest periodic orbit. Our best upper bounds on these averages are all within 1\% of being sharp, which would be nearly impossible without computer assistance.

A main contribution of the present work, aside from producing novel bounds for the Lorenz system, is that we extend the methods of \cite{Chernyshenko2014, Fantuzzi2016} to produce bounds that are mathematically rigorous, including some that depend analytically on the parameter $r$. Bounds implied by numerical SDP solutions are not rigorous because such solutions typically violate their equality and inequality constraints by small margins. These slight inaccuracies can easily suggest false bounds because the numerical conditioning of the SDPs is often poor. Here we employ two complementary methods of proving rigorous bounds using SDPs. The first approach, described and illustrated in \S\ref{sec: non-sharp}, is to compute an approximate numerical solution and then use interval arithmetic to bound how far it is from the true solution. This procedure is automated by the software package VSDP \cite{Jansson2006} and can produce very tight bounds, but it can never produce perfectly sharp bounds because of the inherent conservativeness  of interval arithmetic. The second approach, described and illustrated in \S\ref{sec: sharp}, is to solve the relevant SDPs exactly. In certain situations exact SDP solutions yield sharp bounds. Here we obtain sharp bounds by solving relatively small SDPs analytically with computer assistance.

Section \ref{sec: method} summarizes the formulation of SDPs that imply bounds on time averages and then extends this framework to parameter-dependent bounds. Section \ref{sec: Lorenz} reviews some properties of trajectories of the Lorenz system. Such knowledge of trajectories is not needed to construct bounds and is often unavailable, but since the Lorenz system has been studied extensively we can use this information to judge the quality of our bounds. Sections \ref{sec: non-sharp} and \ref{sec: sharp} describe methods of obtaining rigorous bounds from SDPs using interval arithmetic and exact solutions, respectively, as well as reporting the bounds we have proved for the Lorenz system by each approach. Bounds reported in \S\ref{sec: non-sharp} are not sharp and apply only at the standard chaotic parameter values, while bounds reported in \S\ref{sec: sharp} are sharp and apply over ranges of parameter values. Section \ref{sec: con} summarizes our best bounds for the Lorenz system and discusses the potential of the SDP bounding methodology.

\section{\label{sec: method}Bounding time averages using semidefinite programming}

Consider a well-posed finite-dimensional dynamical system
\beq
\tfrac{d}{dt}\x = \f(\x), \quad \x,\f\in\mathbb R^n
\label{eq: generic ode}
\eeq
that is bounded, meaning all trajectories $\x(t)$ remain bounded as $t\to\infty$. For any function $\varphi(\x)$, let $\ov\varphi$ denote its infinite-time average along a trajectory:
\beq
\ov \varphi = \limsup_{T\to\infty}\frac{1}{T}\int_0^T\varphi(\x(t))dt.
\label{eq: overbar}
\eeq
Our results are unchanged if averages are defined using $\liminf$ instead of $\limsup$. The value of $\ov\varphi$ may depend on the trajectory $\x(t)$ that is averaged over. Here we construct bounds on $\ov\varphi$ that are global, meaning they apply to all trajectories. In other words, we construct lower bounds $L$ and upper bounds $U$ such that
\beq
L \le \min_{\x(t)}\ov\varphi 
\quad \text{ and } \quad
\max_{\x(t)}\ov\varphi \le U.
\label{eq: global bounds}
\eeq
The related challenge of locating the trajectories that extremize $\ov\varphi$ is tackled in \cite{Tobasco2017}. Although we speak of time averages along trajectories, the above extrema are unchanged if they are instead taken over expectations of invariant measures \cite{Tobasco2017}, as studied in the field of ergodic optimization \cite{Jenkinson2006}.

Section \ref{sec: suff} gives sufficient conditions that we use to prove global bounds of the form \eqref{eq: global bounds}. The formulation of SDPs that imply these sufficient conditions is reviewed in \S\ref{sec: relaxation} and extended to parameter-dependent bounds in \S\ref{sec: parameter-dependent}.

\subsection{\label{sec: suff}Sufficient conditions for bounding time averages}

To bound averages over trajectories without knowing the trajectories themselves, we make use of relations between time averages that are implied by the governing system \eqref{eq: generic ode}. Each such relation is satisfied on all trajectories, but none is sufficiently constraining to determine trajectories uniquely. This relaxation is crucial to making the analysis tractable; using the full constraint of the governing system would be tantamount to finding its trajectories exactly, which is generally not possible. Infinitely many constraints can be generated from the fact that time derivatives average to zero in bounded systems: for any differentiable scalar function $V(\x)$ and bounded differentiable trajectory $\x(t)$,
\begin{equation}
\ov{\tfrac{d}{dt}V}= \limsup_{T\to\infty}\frac{1}{T}
	\left[V(\x(T))-V(\x(0))\right] = 0.
\end{equation}
The chain rule gives $\tfrac{d}{dt}V=\f\cdot\nabla V$, so
\beq
\ov{\f\cdot\nabla V}=0.
\label{eq: generic constraint}
\eeq
Different choices for $V$ in the above expression generally provide infinitely many time-averaged relations that hold along every trajectory. However, only particular $V$ yield relations that help to prove a desired bound.

Suppose that for some function $\varphi(\x)$ we want to prove that a lower bound $L\le\ov\varphi$ holds along every trajectory. It would suffice to show that the bound holds pointwise on all trajectories, meaning $L\le\varphi(\x)$ for all $\x\in\mathbb R^n$, but this generally will not be true. Instead we can exploit the identity \eqref{eq: generic constraint} by finding $V$ such that $L\le\varphi(\x)+\f(\x)\cdot\nabla V(\x)$ \emph{does} hold for all $\x\in\mathbb R^n$. This pointwise condition proves the desired results since its time averages is $L\le\ov\varphi$. The pointwise inequality is equivalent to nonnegativity of the function
\begin{equation}
S_L(\x) = \varphi(\x)-L+\f(\x)\cdot\nabla V(\x).
\label{eq: S_L}
\end{equation}
In summary, we have a sufficient condition to prove a lower bound on $\ov\varphi$ for all trajectories:
\beq
\exists V(\x) \text{ such that } S_L(\x) \ge 0 ~\forall \x\in\mathbb R^n
\quad \Longrightarrow \quad
L\le\ov\varphi.
\label{eq: L suff}
\eeq
The lower bound is sharp if $L=\min_{\x(t)}\ov\varphi$.

Analogous arguments for an upper bound lead to a sufficient condition where the function that must be pointwise nonnegative is
\begin{equation}
S_U(\x) = -\left[ \varphi(\x)-U+\f(\x)\cdot\nabla V(\x)\right],
\label{eq: S_U}
\end{equation}
and an upper bound is implied according to
\beq
\exists V(\x) \text{ such that } S_U(\x) \ge 0 ~\forall \x\in\mathbb R^n
\quad \Longrightarrow \quad
\ov\varphi\le U.
\label{eq: U suff}
\eeq
The upper bound is sharp if $U=\max_{\x(t)}\ov\varphi$. The choice of the auxiliary function $V$ is generally different for $S_U$ than for $S_L$.

When applying the sufficient conditions \eqref{eq: L suff} and \eqref{eq: U suff}, it is natural to seek the tightest possible bounds by optimizing over the choice of $V$. It is proved in \cite{Tobasco2017} that such optimization always gives arbitrarily sharp bounds. That is, if $\f$ and $V$ are differentiable on a compact set that all trajectories eventually remain in,
\beq
\sup_{\substack{V(\x)\\S_L(\x)\ge0}}L = \min_{\x(t)}\ov\varphi 
\quad \text{ and } \quad
\max_{\x(t)}\ov\varphi = \inf_{\substack{V(\x)\\S_U(\x)\ge0}}U.
\label{eq: sharpness}
\eeq
The practical challenge is to prove sharp or nearly sharp bounds by constructing auxiliary functions that are as close to optimal as possible.

\paragraph{Simple example.}

Suppose we want to show that $\ov{xy}\ge0$ for all solutions of the Lorenz equations \eqref{eq: Lorenz}. The identity \eqref{eq: generic constraint} with $V=-\tfrac{1}{2\sigma}x^2$ reveals that $\ov{xy}=\ov{x^2}$ on all trajectories. Whereas $xy<0$ on some parts of the attractor, $x^2\ge0$ everywhere in phase space. This implies $\ov{x^2}\ge0$ and thus the desired lower bound $\ov{xy}\ge0$. In the nomenclature of condition \eqref{eq: L suff}, the sufficient condition proving the bound is $S_L=x^2\ge0$.

\subsection{\label{sec: relaxation}Sum-of-squares relaxation and semidefinite programming}

Proving a lower bound using condition \eqref{eq: L suff} entails two related challenges: finding an auxiliary function $V$ that makes $S_L$ nonnegative, and showing that $S_L$ is indeed nonnegative. In this work we consider only polynomial dynamics, meaning that each component of $\f(\x)$ is polynomial in the components of $\x$. We further restrict ourselves to polynomial $V$, so $S_L$ is polynomial also. Even so it is too difficult in general to decide whether $S_L(\x)\ge0$ because the computational complexity of this question is NP-hard \cite{Murty1987}. We thus impose the more tractable condition that $S_L$ can be represented as a sum of squares of other polynomials, which is sufficient for $S_L(\x)\ge0$ but not generally necessary \cite{Hilbert1888}. In turn, the condition that $S_L$ is SOS can be stated as a condition on the matrix representation of $S_L$. This matrix representation, which is called a \emph{Gram matrix}, is any symmetric matrix $\Q$ such that $S_L=\bb^T\Q\bb$, where $\bb(\x)$ is a specified vector of polynomial basis functions. As long as $\bb$ contains enough basis functions, such a $\Q$ always exists and often is not unique. A polynomial is SOS if and only if there exists a Gram matrix that is positive semidefinite, i.e.\ $\Q\succeq0$ (cf.\ \cite{Parrilo2003}). In summary, we replace the sufficient condition of \eqref{eq: L suff} by the stronger condition that $S_L$ is SOS:
\beq
\exists V(\x),\Q,\bb(\x) \text{ such that } S_L=\bb^T\Q\bb \text{ and } \Q\succeq0
\quad \Longrightarrow \quad
L\le\ov\varphi.
\label{eq: L suff sos}
\eeq

In practice the sufficient condition \eqref{eq: L suff sos} is applied by assuming a polynomial ansatz for $V$ with tunable coefficients. This gives an ansatz for $S_L$ in which these coefficients also appear. Based on the monomials in $S_L$ it is simple to choose an adequate polynomial basis vector $\bb$. Once $\bb$ is fixed, the equality $S_L=\bb^T\Q\bb$ furnishes affine constraints on the entries of $\Q$. (Matching coefficients on like monomials gives affine relations involving the entries of $\Q$ and the tunable coefficients of $V$, but the latter can be eliminated if desired.) Thus, for a chosen $V$ ansatz and basis $\bb$, the sufficient condition of \eqref{eq: L suff sos} amounts to affine and semidefinite constraints on $\Q$. Finding a symmetric matrix subject to these two types of constraints is exactly what constitutes an SDP \cite{Boyd2004}. Furthermore, this SDP can be posed as an optimization, searching for the maximum $L$ for which there exists a $\Q$ satisfying condition \eqref{eq: L suff sos}. For a given $V$ ansatz, the best lower bound $L^*\le\ov\varphi$ that can be proved using the SDP formulation is:
\beq
L^* = \displaystyle\max_{\substack{\Q\succeq0\\S_L=\bb^T\Q\bb}} L,
\label{eq: L optimal}
\eeq
where $\bb$ is any fixed basis vector that can represent $S_L$. If the $V$ ansatz is inadequate to prove any finite bound, then $L^*=-\infty$. By analogous arguments, the best upper bound $\ov\varphi\le U^*$ that can be proved with a given $V$ ansatz is
\beq
U^* = \displaystyle\min_{\substack{\Q\succeq0\\S_U=\bb^T\Q\bb}} U.
\label{eq: U optimal}
\eeq
Adding terms to the ansatz of $V$ either improves the optima $L^*$ and $U^*$ or leaves them unchanged, although the approximate optima returned by numerical SDP solvers may change less predictably due to numerical conditioning getting worse as the dimension of $\Q$ increases.

The SDP optimizations \eqref{eq: L optimal} and \eqref{eq: U optimal} are convex: the matrices $\Q$ satisfying given affine and semidefinite constraints form a convex set. Various solvers are available to compute approximate numerical solutions, such as SDPT3 \cite{Tutuncu2003} or Mosek \cite{Mosek}. Numerical approximations to $L^*$ and $U^*$ returned by SDP solvers do not constitute mathematically rigorous bounds because of roundoff error. Methods of obtaining rigorous results by using interval arithmetic and by finding exact SDP solutions are described at the beginnings of \S\ref{sec: non-sharp} and \S\ref{sec: sharp}, respectively. Furthermore, trial-and-error with numerical solutions can simplify subsequent rigorous analyses by suggesting terms that can be omitted from $V$ and $\bb$.

\subsection{\label{sec: parameter-dependent}Parameter-dependent bounds}

When a dynamical system is parameterized by some vector of parameters $\p$,
\beq
\tfrac{d}{dt}\x = \f(\x,\p), \quad \x,\f\in\mathbb R^n, \quad \p\in\mathbb R^m,
\label{eq: param ode}
\eeq
we can seek bounds that depend analytically on parameters. Bounds that are polynomial in the components of $\p$ can be constructed using the SDP framework. In the case of a polynomial lower bound $L(\p)$, the sufficient condition \eqref{eq: L suff sos} is extended by letting the auxiliary function $V$ and basis $\bb$ be polynomials in the components of $\p$ as well as $\x$:
\beq
\begin{array}{r}
\exists V(\x,\p),\Q,\bb(\x,\p) \text{ such that } 
S_L=\bb^T\Q\bb \text{ and } \Q\succeq0 
\quad \Longrightarrow \quad
L(\p)\le\displaystyle\ov\varphi.
\end{array}
\label{eq: L param sos}
\eeq
The bound is sharp if $L(\p)=\min_{\x(t)}\ov\varphi$ for all $\p$. Once a $V$ ansatz and basis $\bb$ are specified, the sufficient condition \eqref{eq: L param sos} is an SDP. This condition can be posed as an optimization by maximizing a coefficient in the lower bounds ansatz.

At least two obstacles can arise when treating parameters analytically, as opposed to fixing them numerically. The first occurs when trying to prove bounds that fail for some parameter values. An example for the Lorenz system is the bound $\ov{z^3}\le(r-1)^3$ proved in \S\ref{sec: z3} that holds when $r\ge1$. Since the bound fails for some $r<1$ it cannot be proved by condition \eqref{eq: L param sos} if $r$ is regarded as an arbitrary variable. In this instance the obstacle can be surmounted by instead proving $(r-1)\ov{z^3}\le(r-1)^4$, which is true for all $r$ and implies the desired result. A second obstacle is that optimal choices of $L(\p)$ or $V(\x,\p)$ may be non-polynomial in $\p$. This occurs in \S\ref{sec: sharp} when we choose $V$ that are polynomial in $r$ but non-polynomial in $\beta$ and $\sigma$. In this case $r$ is the only parameter that can be included analytically in the polynomial basis vector $\bb$, so the Gram matrix $\Q$ generally must depend on $\beta$ and $\sigma$. Parameters on which $\Q$ depends must be fixed numerically, except in SDPs small enough to permit analytical study.

\section{\label{sec: Lorenz}The Lorenz system: trajectories and averages}

Let us recall the phase space of the Lorenz system \eqref{eq: Lorenz}. All trajectories are bounded forward in time, approaching a global attractor whose location can be approximated by finding trapping regions that trajectories do not leave \cite{Doering1995b, Swinnerton-Dyer2001}. Certain trapping regions reveal that $z\ge0$ on the attractor for all positive values of the parameters $(\beta,\sigma,r)$, meaning that as $t\to\infty$ every trajectory approaches a point or set of points in the half-space $z\ge0$.

The Lorenz system has three equilibria: one at $\x=(0,0,0)$ that exists for all parameter values and two nonzero equilibria,
\begin{equation}
\x^\pm=\big(\pm\sqrt{\beta(r-1)},\pm\sqrt{\beta(r-1)},r-1\big),
\label{eq: nonzero}
\end{equation}
that exist when $\beta(r-1)>0$. The zero equilibrium is an example of a symmetric trajectory, meaning it is invariant under the symmetry $(x,y)\mapsto(-x,-y)$ of the governing equations. The nonzero equilibria $\x^\pm$ comprise a pair of asymmetric trajectories that are mapped to one another by this symmetry. The stability of the equilibria and the existence of other invariant structures depend on the parameter values.

Section \ref{sec: Lorenz phase} details the phase space of the Lorenz system at the standard chaotic parameter values. Section \ref{sec: Lorenz moments} discusses the time-averaged moments we bound here and some exact relations between them. Section \ref{sec: extremal traj} reports which trajectories appear to maximize or minimize these mean moments at the standard parameter values -- observations that help us judge the sharpness of the bounds reported in \S\S\ref{sec: non-sharp}--\ref{sec: sharp}.

\subsection{\label{sec: Lorenz phase}Phase space at the standard parameters}

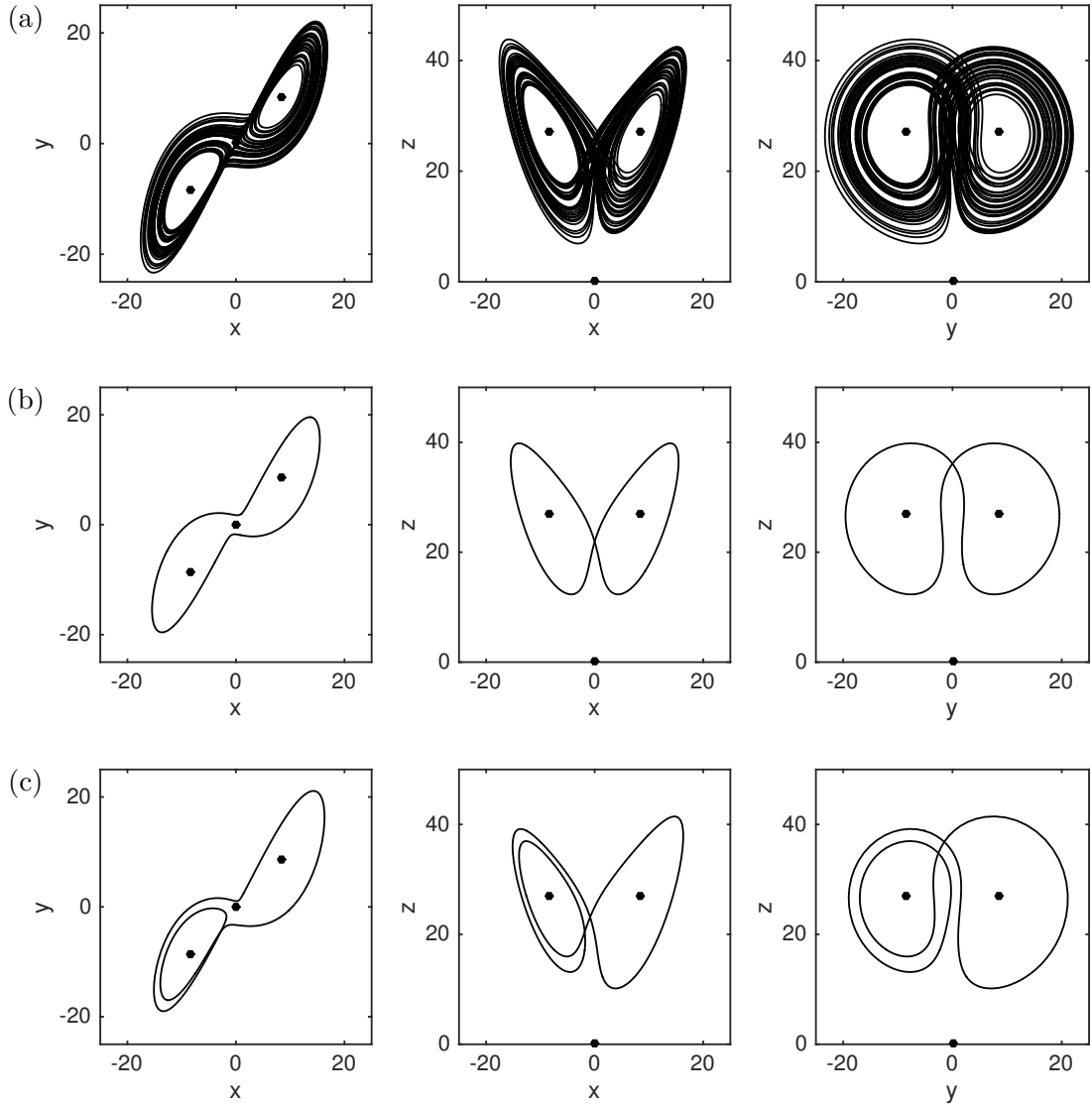
\begin{figure}[p]
\begin{center}
\begin{tikzpicture}
\node at (0,0) {\includegraphics[width=410pt,trim={50pt 40pt 60pt 50pt},clip]{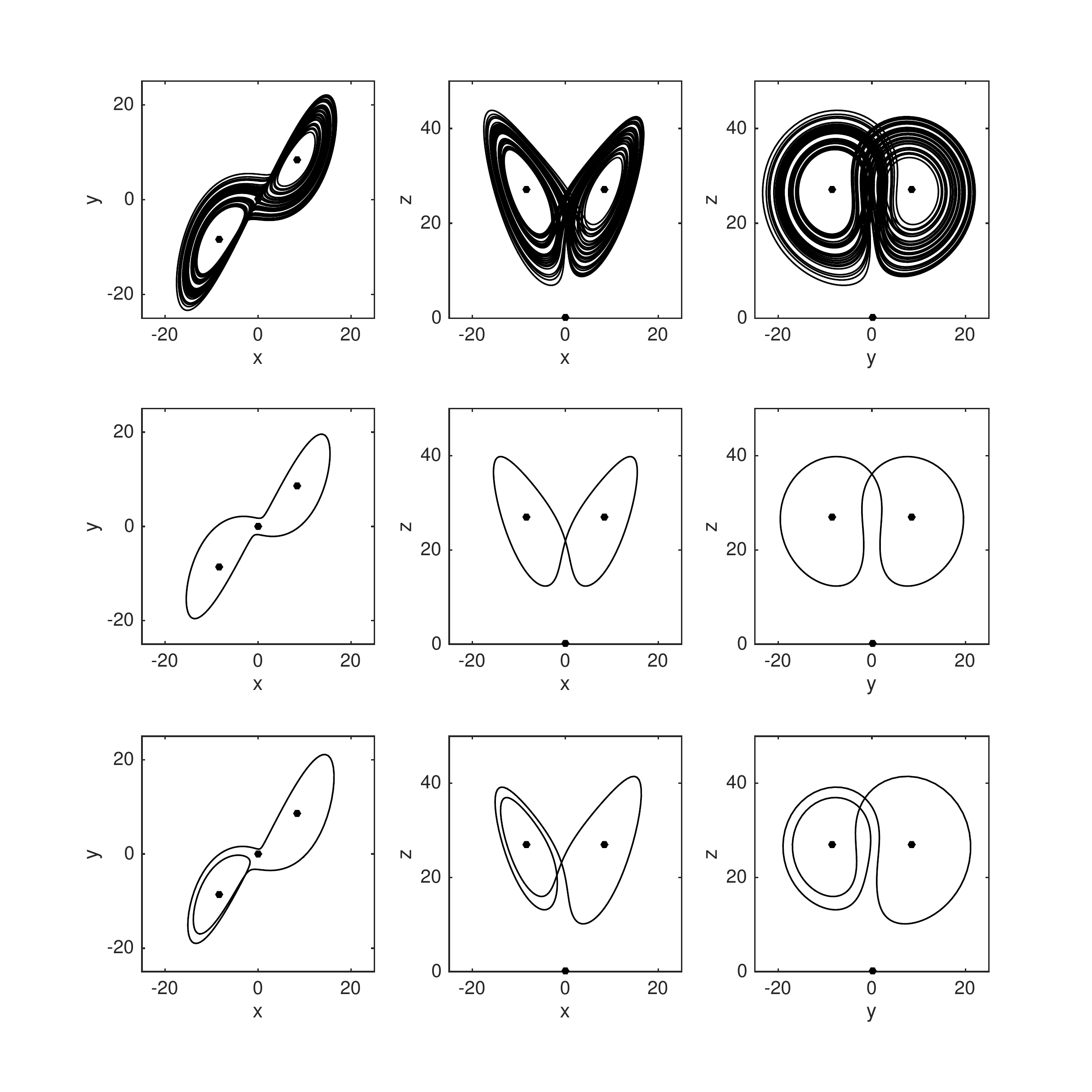}};
\draw (-7.2,7.2) node {(a)};
\draw (-7.2,2.08) node {(b)};
\draw (-7.2,-3.04) node {(c)};
\end{tikzpicture}
\end{center}
\caption{\label{fig: traj} Trajectories of the Lorenz system at the standard parameters $(\beta,\sigma,r)=(8/3,10,28)$: (a) part of a chaotic trajectory, (b) the symmetric periodic orbit $\oplus\ominus$, and (c) the asymmetric periodic orbit $\ominus$$\ominus$$\oplus$. The equilibria at the origin and at $\x^\pm$ are also shown. Periodic orbits are from numerical computations of Viswanath \cite{Viswanath2003}.}
\end{figure}

The phase space of the Lorenz system \eqref{eq: Lorenz} at the chaotic parameters $(\beta,\sigma,r)=(8/3,10,28)$ originally considered by Lorenz \cite{Lorenz1963} is very well studied. All three equilibria are unstable, as are the infinitely many periodic orbits. The set of initial conditions lying on the equilibria or periodic orbits, or on their stable manifolds, has zero volume. Generic initial conditions instead produce chaotic trajectories, quickly converging to the strange attractor whose existence is proven \cite{Tucker1999}. Figure~\ref{fig: traj}a shows a piece of a chaotic trajectory projected onto the three coordinate planes. The complicated strange attractor is built from simpler invariant structures: the zero equilibrium and every periodic orbit are part of the strange attractor, as are their unstable manifolds, and any chaotic trajectory eventually passes arbitrarily close to all of these structures (although close passes to the origin are very rare). The nonzero equilibria $\x^\pm$ are not embedded in the strange attractor but are nearby; the Euclidean distance between $\x^\pm$ and the unstable manifold of the origin is about 1.56.

Each periodic orbit can be labelled by a periodic sequence of $\oplus$ and $\ominus$ symbols, corresponding to the sequence in which it winds around the nonzero equilibria $\x^+$ and $\x^-$, respectively. Winding can be defined precisely using a Poincar\'e section such as the plane $z=r-1$. At the standard parameters, this symbolic dynamics labels each orbit uniquely \cite{Sparrow1982}. For some symbol sequences there does not exist a corresponding orbit, including the one-symbol orbits $\oplus$ and $\ominus$ that would periodically circle $\x^+$ and $\x^-$, respectively. Thus every periodic orbit winds around both $\x^+$ and $\x^-$ at least once. Figure~\ref{fig: traj}b shows the shortest orbit, which has the periodically repeating sequence $\oplus\ominus$ and is symmetric. The second-shortest orbits are asymmetric and have periodic sequences $\oplus$$\oplus$$\ominus$ and $\ominus$$\ominus$$\oplus$. Figure~\ref{fig: traj}c shows the $\ominus$$\ominus$$\oplus$ orbit. (If desired, one can define a symmetry-invariant symbolic dynamics where orbits like $\oplus$$\oplus$$\ominus$ and $\ominus$$\ominus$$\oplus$ are identified \cite{Eckhardt1994}.)

Since the bounds we construct here are global, they apply not only to the generic chaotic trajectories but also to all the non-generic trajectories that tend to equilibria or periodic orbits. For every average quantity $\ov\varphi$ that we have examined, the chaotic value is \emph{not} extremal; there always exist periodic orbits, and sometimes also equilibria, on which $\ov\varphi$ is larger or smaller than its chaotic value. (This is consistent with the interpretation of chaotic averages as weighted averages over periodic orbits \cite{Cvitanovic2016}.) Therefore sharp global bounds give estimates of chaotic averages that are at least slightly conservative.

\subsection{\label{sec: Lorenz moments}Symmetric moments}

In this work we bound averages of moments that are symmetric under $(x,y)\mapsto(-x,-y)$, meaning monomials $x^ly^mz^n$ where $l,m,n\ge0$ are integers and $l+m$ is even. Antisymmetric moments, where $l+m$ is odd, average to zero on chaotic trajectories, although not on asymmetric trajectories such as the equilibria $\x^\pm$ or the periodic orbit $\ominus$$\ominus$$\oplus$ in figure~\ref{fig: traj}c. Here we bound the moments listed in table~\ref{tab: sym mom}, which are all of the symmetric moments whose degree $l+m+n$ is no larger than four.

\begin{table}[t]
\begin{center}
\caption{\label{tab: sym mom}The eighteen moments $\mom$ up to quartic degree that are invariant under the symmetry $(x,y)\mapsto(-x,-y)$ of the Lorenz equations.}
\begin{tabular}{cl}
Degree & Symmetric moments \\ 
\hline
1 & $z$ \\ 
2 & $x^2$, $xy$, $y^2$, $z^2$ \\ 
3 & $x^2z$, $y^2z$, $xyz$, $z^3$ \\ 
4 & $x^4$, $x^3y$, $x^2y^2$, $x^2z^2$, $xy^3$, $xyz^2$, $y^4$, $y^2z^2$, $z^4$
\end{tabular}
\end{center}
\end{table}

Average moments can vary between different trajectories. All moments vanish at the zero equilibrium. At the nonzero equilibria \eqref{eq: nonzero} symmetric moments take the values
\begin{equation}
x^ly^mz^n\big|_{\x^\pm} = \beta^{(l+m)/2}(r-1)^{(l+m)/2+n}.
\label{eq: mom nonzero}
\end{equation}
The value of $\ov\mom$ is the same on every chaotic trajectory at parameters where the strange attractor is ergodic, whereas it generally varies between different periodic orbits.

For convenience we define the normalized mean moment $\nmom$ on a trajectory as the time average $\ov\mom$ on that trajectory, divided by the moment's value \eqref{eq: mom nonzero} at the nonzero equilibria,
\begin{equation}
\nmom =\frac{\ov\mom}{\mom\big|_{\x^\pm}}.
\label{eq: M}
\end{equation}
Bounds are often reported here in terms of $\nmom$, although we constructed them using the unnormalized averages $\ov\mom$. For trajectories on the stable manifolds of the zero or nonzero equilibria, $\nmom =0$ or $\nmom =1$, respectively. Other values of $\nmom$ occurring on chaotic or periodic trajectories cannot be calculated exactly but can be approximated numerically.

Some mean moments of the Lorenz system are exactly proportional and thus have identical normalized values. This can be deduced \emph{a priori} from the general identity \eqref{eq: generic constraint} by choosing $V=\tfrac{1}{n\sigma}x^n$ and $V=\tfrac{1}{n}z^n$, respectively, to find that
\begin{align}
\ov{x^{n-1}y} &= \ov{x^n}, & \ov{xyz^{n-1}} &= \beta\ov{z^n}
\label{eq: prop families}
\end{align}
on each trajectory for all $n\ge1$. The first few equalities in these families ($n=1,2,3$) relate time averages of some of the symmetric moments studied here:
\begin{align}
\ov{x^2} &= \ov{xy} = \beta\ov{z}, &
\ov{x^4} &= \ov{x^3y}, &
\ov{xyz} &= \beta\ov{z^2}, &
\ov{xyz^2} &= \beta\ov{z^3}.
\label{eq: proportional}
\end{align}
Among the above moments, we construct bounds on only $\ov z$, $\ov{z^2}$, $\ov{z^3}$, and $\ov{x^4}$ since these imply bounds on the other five moments also.

Beyond the proportionality relations \eqref{eq: proportional}, there exist infinitely many relations involving three or more mean moments. Some of these underlie our present bounding methods, where we seek relations $\ov\varphi=\ov{\varphi+\f\cdot\nabla V}$ such that $\varphi(\x)+\f(\x)\cdot\nabla V(\x)$ obeys the desired bounds pointwise. Other polynomial relations can be useful in other ways. For instance, the equalities of table~\ref{tab: all relations} in Appendix \ref{app: equalities} reveal that the mean values of all eighteen moments listed in table~\ref{tab: sym mom} are linear combinations of just six: $\ov z$, $\ov{z^2}$, $\ov{y^2z}$, $\ov{z^3}$, $\ov{y^2z^2}$, and $\ov{z^4}$. Such relations are useful for computing exact values of mean moments in terms of other exact values, but they are of limited use in deriving sharp bounds. To see why, consider the relation $\ov{y^2}=\beta\big(r\ov z-\ov{z^2}\big)$. On any particular trajectory, $\ov z$ and $\ov{z^2}$ determine $\ov{y^2}$. However, sharp bounds on $\ov z$ and $\ov{z^2}$ do \emph{not} imply a sharp upper bound on $\ov{y^2}$. This is because the maxima of $\ov{y^2}$, $\ov z$, and $-\ov{z^2}$ occur on different trajectories -- on the $\oplus$$\ominus$ orbit, the $\x^\pm$ equilibria, and the zero equilibrium, respectively. Bounding $\ov{y^2}$ directly using the SDP \eqref{eq: U optimal} gives a better result.

\subsection{\label{sec: extremal traj}Conjectured extremal trajectories}

Before constructing bounds, we can obtain conjectures for the maximum and minimum averages of the moments in table \ref{tab: sym mom} by searching among a large number of candidate trajectories. These trajectories include the three equilibria, a chaotic trajectory, and 1424 different periodic orbits. The chaotic trajectory was obtained by numerically integrating for time $10^7$ using fourth-order Runge--Kutta time steps of size $10^{-3}$; it winds around the nonzero equilibria more than 13 million times. The periodic orbits were computed by Viswanath \cite{Viswanath2003, Viswanath2004}. They include all orbits up to thirteen symbols in length and some longer ones with sequences of the form $\oplus^N\ominus$ or $\oplus^N\ominus^N$.

Among the eighteen symmetric moments considered, eight of them ($z$, $x^2$, $xy$, $z^2$, $x^2z$, $xyz$, $z^3$, $xyz^2$) appear to be maximized at the nonzero equilibria $\x^\pm$, meaning that $\max_{\x(t)} \nmom =1$. These conjectures can be confirmed by sharp upper bounds on $\ov{z}$, $\ov{z^2}$, $\ov{z^3}$, and $\ov{x^2z}$ since the other four moments are proportional to these according to \eqref{eq: proportional}. Sharp bounds indeed have been proved for $\ov{z}$ by Malkus \cite{Malkus1972}, for $\ov{z^2}$ by Knobloch \cite{Knobloch1979}, and for $\ov{z^3}$ in \S\ref{sec: z3}. A sharp upper bound has not been proved for $\ov{x^2z}$, although the bound reported in \S\ref{sec: enclosures} is extremely close. The other ten moments listed in table \ref{tab: sym mom} ($y^2$, $y^2z$, $x^4$, $x^3y$, $x^2y^2$, $x^2z^2$, $xy^3$, $y^4$, $y^2z^2$, $z^4$) all appear to be maximized by the $\oplus\ominus$ periodic orbit of figure~\ref{fig: traj}b. In these cases $\max_{\x(t)} \nmom >1$. (In fact we find that all symmetric moments of degrees five through nine are also maximized on the $\oplus\ominus$ orbit, while some higher-degree moments including $x^2y^8$, $xy^9$, and $y^{10}$ are maximized on longer periodic orbits.) For each of these ten moments we report an upper bound in \S\ref{sec: non-sharp upper} that is within 1\% of the corresponding average on the $\oplus\ominus$ orbit. The minima of all eighteen mean moments are attained at the origin, meaning that $\min_{\x(t)} \nmom =0$, as confirmed by the sharp lower bounds of \S\ref{sec: lower}.

\section{\label{sec: non-sharp}Non-sharp bounds using interval arithmetic}

One way to construct a rigorous bound using SDP optimization is to find an approximate numerical solution and then apply perturbation methods, made rigorous by interval arithmetic, to estimate how far the true optimum is from the numerical one. Here we do so using the software VSDP \cite{Jansson2006}, which automates this approach. Section \ref{sec: rigorous pert} describes the method, including why the resulting bounds can be very tight but never perfectly sharp. Section \ref{sec: non-sharp Lorenz} reports upper bounds found by applying this method to the Lorenz system at the standard parameters. Section \ref{sec: why} explains why tight bounds are easier to prove on some quantities than on others, both in the Lorenz example and in general.

\subsection{\label{sec: rigorous pert}Perturbation methods with interval arithmetic}

When using the SOS framework to prove an upper bound, for instance, the best bound produced by a given ansatz of the auxiliary function $V$ is the exact optimum $U^*$ of the SDP in \eqref{eq: U optimal}. Numerical SDP solvers only approximate this optimum and the matrix $\Q^*$ that achieves it. However, perturbation methods can be applied to this approximate $\Q^*$ to find an enclosure $[U^-,U^+]$ that contains the true optimum $U^*$. The software VSDP implements this procedure by employing an external SDP solver (we use SDPT3~\cite{Tutuncu2003}) and the interval arithmetic package INTLAB \cite{Rump1999} to rigorously enclose $U^*$. When VSDP succeeds, the result is a verified upper bound $\ov\varphi\le U^+$. This approach works well in a number of situations but has two main limitations. The first is purely numerical: if the dimension of phase space and/or the polynomial degree of $S_U$ are too large, VSDP may return a value of $U^+$ that is much larger than $U^*$ or infinite. The second limitation is that $U^+$ is always strictly larger than $U^*$. This is not important in the many cases where $U^+$ is close to $U^*$, and $U^*$ is not a sharp bound anyway. In cases where the exact $U^*$ gives a sharp bound, however, it is desirable to prove it by solving the SDP exactly. Such sharp bounds are the topic of \S\ref{sec: sharp}.

\subsection{\label{sec: non-sharp Lorenz}Non-sharp bounds for the Lorenz system at the standard parameters}

For the Lorenz system at the standard chaotic parameter values we have constructed upper bounds on all eighteen moments listed in table~\ref{tab: sym mom} by applying the software VSDP to the SDP \eqref{eq: U optimal}. Results are reported as bounds on the normalized mean moments $\nmom$ defined by \eqref{eq: M}. We have constructed bounds using $V$ ans\"atze of degree 2, 4, 6, 8, and 10. Each ansatz for $V$ includes all symmetric monomials except some that can be excluded at the highest degree (cf.\ \S\ref{sec: Lorenz formulation} below).
 
Recall from \S\ref{sec: extremal traj} that among the eighteen mean moments we consider, eight appear to be maximized on the nonzero equilibria $\x^\pm$, while the other ten appear to be maximized on the $\oplus$$\ominus$ periodic orbit. For the first eight, where it seems that $\max_{\x(t)}\nmom =1$, table \ref{tab: enclosures} reports enclosures $\left[U^-,U^+\right]$ that contain upper bounds on $\nmom$. For the other ten moments, table \ref{tab: non-sharp} reports verified upper bounds $U^+$, along with the maximum known averages. The results in tables \ref{tab: enclosures} and \ref{tab: non-sharp} are discussed in \S\ref{sec: enclosures} and \S\ref{sec: non-sharp upper}, respectively. Appendix \ref{app: vsdp} gives computational details.

\begin{table}[t]
\begin{center}
\caption{\label{tab: enclosures}Enclosures of upper bounds on normalized mean moments $\nmom$ of the Lorenz system at the standard parameter values $(\beta,\sigma,r)=(8/3,10,28)$, computed using $V$ ans\"atze of degree 2, 4, 6, and 8. Each tabulated moment appears to be maximized at the nonzero equilibria $\x^\pm$, in which case the sharp upper bound would be $\nmom\le1$. The verified upper bounds are the $U^+$ values. Underlined digits indicate agreement between $U^-$ and~$U^+$. Moments listed together have the same normalized means according to \eqref{eq: proportional}.}
\begin{tabular}{lp{75pt}p{75pt}p{75pt}p{75pt}}
Moment & \multicolumn{4}{c}{Enclosure {\small $\begin{bmatrix}U^-\\U^+\end{bmatrix}$} containing upper bound on $\nmom$} \\[9pt]
\hhline{~----}
& degree 2 & degree 4 & degree 6 & degree 8  \\
\hline\\[-10pt]
$z$, $x^2$, $xy$ 	
	& $\begin{bmatrix}\underline{0.999999999}89\\
				    \underline{1.000000000}09\end{bmatrix}$ \\
$z^2$, $xyz$	
	& $\begin{bmatrix}\underline{0.9999999999}4\\
				    \underline{1.0000000000}4\end{bmatrix}$ \\
$z^3$, $xyz^2$ 	
	&& $\begin{bmatrix}\underline{0.99999999}86\\
				       \underline{1.00000000}02\end{bmatrix}$ \\
$x^2z$
	&& $\begin{bmatrix}\underline{1.00236685}1\\
				       \underline{1.00236685}3\end{bmatrix}$
	& $\begin{bmatrix}\underline{1.0006603}2\\
				    \underline{1.0006603}9\end{bmatrix}$
	& $\begin{bmatrix}\underline{0.999}0309\\
				    \underline{1.000}0003\end{bmatrix}$ 
\end{tabular}
\end{center}
\end{table}

\begin{table}[t]
\begin{center}
\caption{\label{tab: non-sharp}Verified upper bounds $U^+$ on normalized mean moments $\nmom$ of the Lorenz system at the standard parameter values $(\beta,\sigma,r)=(8/3,10,28)$. The tabulated moments are all symmetric moments up to quartic degree that are \emph{not} maximized by the nonzero equilibria $\x^\pm$. Bounds are shown for $V$ ans\"atze of degrees 2 through 10, along with the maximum known values of $\nmom$, all of which occur on the $\oplus$$\ominus$ periodic orbit. Underlined digits of bounds agree with these maximum known values. The moments $x^4$ and $x^3y$ have the same normalized means according to \eqref{eq: proportional}.}
\begin{tabular}{lp{50pt}p{50pt}p{50pt}p{50pt}p{50pt}l}
Moment & \multicolumn{5}{c}{Upper bound} & Maximum \\
\hhline{~-----~}
& degree 2 & degree 4 & degree 6 & degree 8	& degree 10 \\
\hline \\[-12pt]
$y^2$ 				& 7.2593
					& \underline{1.2}585
					& \underline{1.16}94
					& \underline{1.162}7
					& \underline{1.16}49 
					& 1.1621684 \\
$y^2z$ 			&& \underline{1.0}480
					& \underline{1.040}4
					& \underline{1.039}6
					& \underline{1.039}7 
					& 1.0394975 \\
$x^4$, $x^3y$	&& \underline{2}.5702
					& \underline{2}.1334	
					& \underline{1.9}318	
					& \underline{1.91}64 
					& 1.9111906 \\ 
$x^2y^2$			&& 3.8772
					& \underline{2}.7756	
					& \underline{2.3}514	
					& \underline{2.3}220 
					& 2.2975630 \\ 
$x^2z^2$			&& \underline{1}.2822	
					& \underline{1}.2053
					& \underline{1.19}05
					& \underline{1.189}9 
					& 1.1893425 \\ 
$xy^3$  			&& \underline{4}.7666
					& \underline{3}.9332
					& \underline{3}.1236
					& \underline{3.0}239 
					& 2.9987454 \\ 
$y^4$ 				&& 18.766
					& 6.1518
					& \underline{4}.4757
					& \underline{4.1}842 
					& 4.1459937 \\ 
$y^2z^2$			&& \underline{1}.1226
					& \underline{1.0}640
					& \underline{1.049}2
					& \underline{1.048}9 
					& 1.0484088 \\
$z^4$				&& \underline{1.1}966
					& \underline{1.11}99
					& \underline{1.115}8
					& \underline{1.11}68 
					& 1.1155092 
\end{tabular}
\end{center}
\end{table}

\subsubsection{\label{sec: enclosures}Moments maximized by the nonzero equilibria}

For the mean moments that are maximized on the nonzero equilibria $\x^\pm$, it is possible to prove the sharp upper bounds $\nmom \le1$ using the SDP \eqref{eq: U optimal} with $V$ of finite degree. Any upper enclosure $U^+$ verified by VSDP will be strictly larger than $U^*$ and thus not a sharp bound. Nonetheless, the enclosures shown in table~\ref{tab: enclosures} give verified bounds $U^+$ that are very close to sharp, and they inform analytical proofs by suggesting what terms are needed in $V$. These enclosures suggest that the exact optima $U^*$ yield sharp bounds on $\ov{z}$ and $\ov{z^2}$ for $V$ of degree 2, on $\ov{z^3}$ for $V$ of degree 4, and on $\ov{x^2z}$ for $V$ of degree 8. At these degrees, each enclosure is very narrow and includes the value $U=1$ that would be sharp. The apparent sharpness of the exact $U^*$ is confirmed for $\ov z$ and $\ov{z^2}$ by past analytical results \cite{Malkus1972, Knobloch1979} and for $\ov{z^3}$ by the bound we prove in \S\ref{sec: z3} using quartic $V$. For the mean moment $\ov{x^2z}$, it is evident that $V$ of degree 4 or 6 are insufficient to prove the sharp bound since $U^->1$. (On the other hand, at some different values of $\beta$ and $\sigma$ we have found that quartic $V$ gives sharp bounds on $\ov{x^2z}$.) We have not proved analytically that $V$ of degree 8 indeed yields $U^*=1$, so we settle for the slightly imperfect bound $\overline{\raisebox{0pt}[\height-.8pt]{$\overline{x^2z}$}}\le1.0000003$ verified by VSDP.

\subsubsection{\label{sec: non-sharp upper}Moments maximized by the $\oplus$$\ominus$ periodic orbit}

For the mean moments in table~\ref{tab: non-sharp}, all of which are apparently maximized by the $\oplus$$\ominus$ periodic orbit, the SDP approach cannot produce a perfectly sharp upper bound for reasons explained in the next subsection. At best, the optimum $U^*$ of the SDP \eqref{eq: U optimal} may approach a sharp bound as the degree of $V$ approaches infinity. The verified bound $U^+$ returned by VSDP does not become sharp in this limit since enclosures become more conservative or infinite as SDPs grow in size. In practice the smallest value of $U^+$ is achieved by a $V$ ansatz of intermediate degree; the verified bounds in table~\ref{tab: non-sharp} continue to improve up to degree 8 or 10. For each tabulated moment, the best bound is within 1\% of being sharp.

\subsection{\label{sec: why}Constraints on $S_U$ along trajectories}

To understand why some quantities are harder to bound than others, suppose one wants to prove an upper bound $\ov\varphi\le U$ using the sufficient condition \eqref{eq: U suff}. The function $S_U(
\x)$ must be nonnegative everywhere, but it also cannot be too large; on each trajectory its average is the margin between $\ov\varphi$ and the bound:
\beq
\ov{S_U} = U-\ov\varphi.
\label{eq: S_U mean}
\eeq
The closer a trajectory comes to saturating the bound $U$, the more strongly $S_U$ is constrained along it. (This observation underlies a method of using the level sets of $S_U$ to locate the trajectories that maximize $\ov\varphi$ \cite{Tobasco2017}.)

If $U$ is a sharp bound, then $S_U$ must vanish on every invariant trajectory that saturates the bound, as follows from $\ov{S_U}=0$ and $S_U(\x)\ge0$. This is why our SOS methods, when applied to the Lorenz system, do not give sharp bounds on averages maximized by the $\oplus$$\ominus$ orbit. We expect that the $\oplus$$\ominus$ orbit is a non-algebraic curve, in which case it is impossible for the polynomial $S_U$ to vanish everywhere along it, as observed in \cite{Fantuzzi2016}. On the other hand, SOS methods can produce sharp bounds on averages that are extremized by equilibria, including sharp upper bounds in the Lorenz system on quantities maximized at $\x^\pm$. The requirement that $S_U$ (and $\nabla S_U$) vanish at a finite number of maximizing equilibria can be satisfied by polynomial $S_U$. Often what demands higher-degree $V$ is the constraints on $S_U$ not at these equilibria but elsewhere in phase space. One indicator of how strongly $S_U$ is constrained elsewhere is the largest value of $\ov\varphi$ -- that is, the smallest value of $\ov{S_U}$ -- among periodic orbits. For $\ov{z^3}$ and $\ov{x^2z}$ in the Lorenz example, the largest such values we have found occur on the $\oplus^{12}$$\ominus^{12}$ and $\oplus$$\ominus$ orbits, respectively, and are smaller than the maxima at $\x^\pm$ by 5.6\% and 3.3\%. The constraint on $S_U$ is stronger in the latter case, and indeed table~\ref{tab: enclosures} suggests that sharp bounds on $\ov{z^3}$ and $\ov{x^2z}$ require $V$ of degree 4 and 8, respectively.

\section{\label{sec: sharp}Sharp bounds using exactly optimal SDP solutions}

Whereas rigorous bounds in \S\ref{sec: non-sharp} are constructed from approximate SDP solutions, the present section describes a different approach in which SDPs are solved exactly. Unlike the approach of \S\ref{sec: non-sharp}, exact solutions can verify not only suboptimal SDP solutions but also exact optima. This is particularly valuable when the optima give sharp bounds, which may be possible when the trajectories saturating the bounds are algebraic sets (cf.\ \S\ref{sec: why}). In the example of the Lorenz system, we can hope to prove sharp lower bounds that are saturated by the zero equilibrium and sharp upper bounds that are saturated by the nonzero equilibria $\x^\pm$. The moments that appear to be maximized on $\x^\pm$, at least at the standard parameters, are the eight listed in table~\ref{tab: enclosures} above.

Section \ref{sec: rigorous exact} discusses how to find exact solutions to SDP optimizations. We then apply this approach to the Lorenz system. Section \ref{sec: Lorenz formulation} explains how the SDP framework can be tailored to exploit some particular features of the Lorenz equations. Sharp upper bounds on $\ov z$ and $\ov{z^2}$ holding for all positive parameters have been proved by previous authors, and \S\ref{sec: z and z2} restates the proof for $\ov{z^2}$ in the SDP framework. The sharp upper bound $\ov{z^3}\le(r-1)^3$ is proved in \S\ref{sec: z3} for all $r\ge1$ and a range of $(\beta,\sigma)$ including the standard values $(8/3,10)$. Sharp lower bounds of zero are constructed in \S\ref{sec: lower}.

\subsection{\label{sec: rigorous exact}Finding exact solutions to SDP optimizations}

Suppose we want to prove a sharp upper bound by verifying the exact optimum $U^*$ of the SDP \eqref{eq: U optimal}. We must find a matrix $\Q^*\succeq0$ that exactly satisfies the relation $S_U=\bb^T\Q^*\bb$ when $U=U^*$. A potential difficulty is that, unless $\bb$ is chosen in a particular way, optimal solutions tend to be only marginally feasible, meaning there exist $\Q^*\succeq0$ that are singular but none that are strictly positive definite. If $\bb$ contains the minimum number of terms needed to represent $S_U$ for any $U$, and if the value of $U$ is fixed rather than optimized, the SDP becomes a feasibility problem that is strictly feasible when $U>U^*$, marginally feasible when $U=U^*$, and infeasible when $U<U^*$. In other words, the matrix $\Q$ can be nonsingular when $U>U^*$ but has additional null vectors at the optimum \cite{Boyd2004}.

Although marginally feasible SDPs sometimes can be solved in practice, it may be easier to replace them with smaller SDPs that are strictly feasible. This requires finding a smaller basis vector $\bhat$ such that when $U=U^*$ there exists $\Qhat\succ0$ satisfying $S_U=\bhat^T\Qhat\bhat$. The reduced basis $\bhat$ is useful for verifying the optimum $U^*$, but the larger basis $\bb$ is still needed when optimizing $U$ since $\bhat$ is insufficient to represent $S_U$ when $U\neq U^*$. Sometimes $\bhat$ can be chosen by inspection, as in \S\ref{sec: z and z2} and \S\ref{sec: z3} below, and otherwise it can be sought with computer assistance \cite{Monniaux2011}.

\subsubsection{Analytical exact solutions}

Very small SDPs can be solved analytically by hand \cite{Powers1998} or using computer algebra \cite{SafeyElDin2010, Henrion2016a}. This approach quickly becomes intractable as $\Q$ grows because the algebraic degree of solutions becomes very large \cite{Nie2010}. (We have had trouble with $8\times8$ matrices.) One advantage of this approach is that $\Q$ can depend analytically on parameters, which is more flexible than having a parameter-dependent basis $\bb$ because the parameter-dependence of $\Q$ does not need to be polynomial. The bounds for the Lorenz system proved in \S\ref{sec: z and z2} and \S\ref{sec: z3} make use of $\Q$ with non-polynomial dependence on $\beta$ and $\sigma$. A second advantage of the analytical approach is that it is not substantially changed if the SDP is only marginally feasible. Even so, we find it convenient in our proofs to replace marginally feasible SDPs with strictly feasible ones, as described in the preceding paragraph.

\subsubsection{Exact solutions in terms of rational numbers}

Exact solutions to SDPs of moderate or large size can be found by a symbolic-numerical approach wherein approximate numerical solutions are projected onto the rational numbers in a way that exactly satisfies the affine and semidefinite constraints. Any parameters on which $\Q$ depends must be fixed, so a bound can be proved only at these values. Rational solutions exist for all strictly feasible SDPs \cite{Peyrl2008}, so in principle any suboptimal bound can be verified. Often there exist optimal solutions $\Q^*$ that are rational also. Section \S\ref{sec: z3 explicit} gives an example of a rational $\Q^*$ that certifies the sharp bound $\ov{z^3}\le(r-1)^3$ in the Lorenz system at the standard parameters $(\beta,\sigma)=(8/3,10)$. In this instance the SDP is small enough that a rational $\Q^*$ can be found analytically. Larger SDPs must be solved numerically and then projected onto the rationals with computer assistance \cite{Peyrl2008, Kaltofen2008, Kaltofen2012, Wu2013}, but we do not implement such methods here.

\subsection{\label{sec: Lorenz formulation}Exploiting the structure of the Lorenz system}

The bounding SDPs \eqref{eq: L optimal} and \eqref{eq: U optimal} often can be tailored to exploit particular features of the equations being studied. In the example of the Lorenz system, we can exploit both the symmetry $(x,y)\mapsto(-x,-y)$ and the nature of the quadratic nonlinearity. The symmetry is shared by the equations and the quantities $\varphi$ being bounded, so we choose ans\"atze of the auxiliary function $V$ to be symmetric also; there is no advantage to a non-symmetric ansatz because the optimal choice of its coefficients would make it symmetric. The functions $S_L$ and $S_U$ defined by \eqref{eq: S_L} and \eqref{eq: S_U} are symmetric also. These polynomials must be nonnegative and thus of even degree. We choose $V$ ans\"atze of even degree. For $S_L$ and $S_U$ to be even also, the quadratic terms of $\f$ must cancel in $\f\cdot\nabla V$, which occurs only when the highest-degree terms of $V$ all take the form $x^p(y^2+z^2)^q$ \cite{Swinnerton-Dyer2001}. Given this constraint and the symmetry condition, the most general bases we can choose for $V$ are
\begin{align}
\text{degree 2: } 
& \left\{ z,x^2,y^2+z^2 \right\} 
\label{eq: deg 2 basis}, \\
\text{degree 4: } 
& \left\{ z, x^2, xy, y^2, z^2, x^2z, xyz, y^2z, z^3, x^4, x^2(y^2+z^2), (y^2+z^2)^2 \right\},
\label{eq: deg 4 basis}
\end{align}
and so on for higher degrees. When $r$ is included as an analytical variable, the highest-degree terms in $V$ must instead take the form $r^sx^p(y^2+z^2-2rz)^q$. When $V$ is quadratic, for instance, the most general basis is $\left\{ z,x^2,y^2+z^2-2rz \right\}$.

Symmetries of a polynomial constrain its Gram matrix to have block diagonal structure, provided the basis elements in $\bb$ are ordered suitably \cite{Gatermann2004, Parrilo2005, Lofberg2009}. In our cases the polynomial $\bb^T\Q\bb$ must be invariant under $(x,y)\mapsto(-x,-y)$. The vector $\bb$ generally contains both symmetric and antisymmetric elements, but for $\bb^T\Q\bb$ to be symmetric there can be no cross-multiplication between the two types. Thus we can split the vector $\bb$ into vectors of symmetric and antisymmetric elements, $\bb_s$ and $\bb_a$, in which case $\Q$ is block diagonal:
\beq
\bb^T\Q\bb
=\begin{bmatrix}\bb_s^T & \bb_a^T\end{bmatrix}
\begin{bmatrix} \Q_s & {\bf 0} \\ {\bf 0} & \Q_a \end{bmatrix}
\begin{bmatrix}\bb_s \\ \bb_a \end{bmatrix}
=\bb_s^T\Q_s\bb_s + \bb_a^T\Q_a\bb_a.
\eeq
The constraint $\Q\succeq0$ then decouples into the computationally easier constraints $\Q_s\succeq0$ and $\Q_a\succeq0$, and the bounding SDPs \eqref{eq: L optimal} and \eqref{eq: U optimal} become
\begin{align}
&
\begin{array}{rl}
L^* = \displaystyle\max_{\Q_s,\Q_a} L 
\quad \text{ subject to } \qquad
S_L=&\hspace{-7pt}\bb_s^T\Q_s\bb_s + \bb_a^T\Q_a\bb_a, \\
\Q_s,\Q_a\succeq&\hspace{-7pt}0,
\end{array}
\label{eq: L optimal Lorenz}
\\[10pt]
&
\begin{array}{rl}
U^* = \displaystyle\min_{\Q_s,\Q_a} U
\quad \text{ subject to } \qquad
S_U=&\hspace{-7pt}\bb_s^T\Q_s\bb_s + \bb_a^T\Q_a\bb_a, \\
\Q_s,\Q_a\succeq&\hspace{-7pt}0.
\end{array}
\label{eq: U optimal Lorenz}
\end{align}
We solve \eqref{eq: L optimal Lorenz} and \eqref{eq: U optimal Lorenz} in the following subsections, along with a version of \eqref{eq: U optimal Lorenz} where $\bb$ depends analytically on the parameter $r$ as described in \S\ref{sec: parameter-dependent} above.

\subsection{\label{sec: z and z2} Analytical upper bounds on $\ov{z}$ and $\ov{z^2}$ in the Lorenz system}

For all $\beta>0$ and $r>1$, the averages $\ov z$ and $\ov{z^2}$ are maximized by trajectories on or tending to the nonzero equilibria $\x^\pm$, meaning that $\ov z\le r-1$ and $\ov{z^2}\le(r-1)^2$. Both bounds can be proved directly in similar ways, as done by Malkus \cite{Malkus1972} and Knobloch \cite{Knobloch1979}, respectively. However, H\"older's inequality and the fact that $z\ge0$ on the global attractor give
\beq
\ov z \le \ov{z^2}^{1/2} \le \ov{z^3}^{1/3}.
\eeq
The bound $\ov z\le r-1$ is thus implied by the bound $\ov{z^2}\le(r-1)^2$, and both would be implied by $\ov{z^3}\le(r-1)^3$. The latter bound on $\ov{z^3}$ is proved in \S\ref{sec: z3}, but only for a subset of positive $(\beta,\sigma)$. In the present subsection we show that $\ov{z^2}\le(r-1)^2$ for all $\beta>0$. The proof is essentially that of Knobloch, but we illustrate how it can be constructed more systematically using the SDP framework.

\subsubsection{Knobloch's proof}

The bound $\ov{z^2}\le(r-1)^2$ can be proved using the sufficient condition \eqref{eq: U suff} with $\varphi=z^2$, $U=(r-1)^2$, and
\beq
V(x,y,z,r) = \tfrac{1}{\beta}\left[ 2z - 2rz + \tfrac{1}{\sigma}x^2+y^2+z^2 \right].
\label{eq: z2 V}
\eeq
These choices define $S_U$ through expression \eqref{eq: S_U}, and the desired bound follows if $S_U\ge0$ for all $(x,y,z,r)$. Finding that $S_U$ has the SOS representation \cite{Knobloch1979}
\beq
S_U(x,y,z,r) = [z-(r-1)]^2 + \tfrac{2}{\beta}(x-y)^2
\label{eq: z2 S_U sos}
\eeq
thus proves $\ov{z^2}\le(r-1)^2$ for all $\beta>0$.

\subsubsection{Systematic construction}

The preceding argument proves the desired bound but does not illustrate how to come up with the choice \eqref{eq: z2 V} for $V$ that makes $S_U$ an SOS polynomial, nor how to find the SOS representation \eqref{eq: z2 S_U sos}. The SDP formulation offers a systematic approach that is needed when $V$ and $S_U$ are more complicated, as in the next subsection. The approach consists of two steps: first solving an SDP optimization numerically to find a $V$ ansatz that seems to give a sharp bound, and then finding the optimal solution analytically.

The SDP we solve in the numerical step is a modification of \eqref{eq: U optimal Lorenz} where the parameter $r$ is included analytically in the basis $\bb$ as described in \S\ref{sec: parameter-dependent}. Anticipating that the sharp bound is $U=(r-1)^2$, we let $U=(r-1)^2+u_0$ and minimize $u_0$:
\begin{align}
\begin{array}{rl}
u^*_0=\displaystyle\min_{\Q_s,\Q_a} u_0
\quad \text{ subject to } \qquad
S_U=&\hspace{-7pt}\bb_s^T\Q_s\bb_s + \bb_a^T\Q_a\bb_a, \\
\Q_s,\Q_a\succeq&\hspace{-7pt}0.
\end{array}
\label{eq: z2 SDP}
\end{align}
The sharp bound can be proved using any $V$ ansatz that gives $u_0^*=0$. The simplest $V$ ansatz to try is one that includes all possible terms up to quadratic degree. As explained above in \S\ref{sec: Lorenz formulation}, the most general choice we must consider is
\beq
V(x,y,z,r) =  c_1 z + c_2x^2 + c_3(y^2+z^2-2rz).
\label{eq z2 V ansatz}
\eeq
Numerically solving the SDP \eqref{eq: z2 SDP} with this $V$ ansatz at the standard values of $(\beta,\sigma)$ gives $u_0^*\approx0$. Surmising that the true optimum is $u_0^*=0$, we move to the analytical step.

With $\varphi=z^2$, $U=(r-1)^2$, and the $V$ ansatz \eqref{eq z2 V ansatz}, the polynomial $S_U$ defined by \eqref{eq: S_U} is
\beq
S_U(x,y,z,r) = (r-1)^2 + c_1\beta z - 2\beta c_3 rz +  2\sigma c_2x^2 
- (c_1 + 2c_2\sigma)xy + 2c_3y^2 + (2\beta c_3-1)z^2.
\label{eq: z2 S_U}
\eeq
This $S_U$ can be represented as $S_U=\bb_s\Q_s\bb_s+\bb_a\Q_a\bb_a$ using the monomial basis vectors
\begin{align}
\bb_s&=\begin{bmatrix}1&r&z\end{bmatrix}^T,
&
\bb_a&= \begin{bmatrix}x&y\end{bmatrix}^T.
\label{eq: z2 monomial}
\end{align}
These bases indeed lead to a successful proof, but our eventual choices of the Gram matrices $\Q_s$ and $\Q_a$ would be singular. As described in \S\ref{sec: rigorous exact} above, it is possible to choose smaller bases $\bhat_s$ and $\bhat_a$ that lead to a simpler calculation and non-singular Gram matrices. A good choice for $\bhat_s$ and $\bhat_a$ follows from the observation in \S\ref{sec: why} that any $S_U$ proving a sharp bound must vanish at the equilibria $\x^\pm$ saturating the bound. This requires that $\bb_s^T\Q_s\bb_s$ vanishes when $z=r-1$ and $\bb_a^T\Q_a\bb_a$ vanishes when $x=y$, so if the optimal $S_U$ can be represented using the vectors \eqref{eq: z2 monomial}, it also can be represented using the single-entry vectors $\bhat_a=[z-(r-1)]$ and $\bhat_a=[x-y]$. That is,
\beq
S_U(x,y,z,r) = \Qhat_s[z-(r-1)]^2 + \Qhat_a(x-y)^2,
\label{eq: z2 S_U 2}
\eeq
where the Gram matrices $\Qhat_s$ and $\Qhat_a$ are scalars in this case.

Equating expressions \eqref{eq: z2 S_U} and \eqref{eq: z2 S_U 2} for $S_U$ gives affine constraints that relate the coefficients $c_i$ of $V$ to $\Qhat_s$ and $\Qhat_a$. These constraints uniquely determine all five values as
\begin{align}
c_1&=2/\beta, &
c_2&=1/\beta\sigma, &
c_3&=1/\beta, &
\Qhat_s&=1, &
\Qhat_a&=2/\beta.
\label{eq: z2 c_i}
\end{align}
In other cases the affine constraints may not fix all coefficients. Although $\Qhat_s$ and $\Qhat_a$ in the present example are determined by the affine constraints alone, they are indeed nonnegative when $\beta>0$, in which case the bound $\ov{z^2}\le(r-1)^2$ is proved. The values given by \eqref{eq: z2 c_i} yield the same choice \eqref{eq: z2 V} for $V$ and SOS representation \eqref{eq: z2 S_U sos} for $S_U$ that underlie Knobloch's proof \cite{Knobloch1979}.

\subsection{\label{sec: z3} Analytical upper bound on $\ov{z^3}$ in the Lorenz system}

In this subsection we prove that $\ov{z^3}$ is maximized by trajectories on or tending to the nonzero equilibria $\x^\pm$, meaning that $\ov{z^3}\le(r-1)^3$, for all $r\ge1$ and a subset of positive $(\beta,\sigma)$ that includes the standard values $(8/3,10)$. For the $(\beta,\sigma)$ where our proof is valid, this result is stronger than the bound $\ov{z^2}\le(r-1)^2$ constructed in the preceding subsection. We treat $r$ analytically as described in \S\ref{sec: parameter-dependent} and thus cannot prove $\ov{z^3}\le(r-1)^3$ directly since the inequality is false for some $r<1$. Instead we prove that $(r-1)\ov{z^3}\le(r-1)^4$ for all $r$, which implies the desired bound on $\ov{z^3}$ when $r\ge1$. We let $\rho=r-1$ and regard $V$ and $S_U$ as polynomials in $(x,y,z,\rho)$, which leads to a simpler proof than using $(x,y,z,r)$. In terms of $\rho$, the bound we construct is $\rho \ov{z^3}\le\rho^4$.

\subsubsection{Numerical determination of an ansatz for $V$}

First we use numerical SDP optimization to find an ansatz for the auxiliary function $V$ that appears to yield a sharp bound at the standard values of $(\beta,\sigma)$, after which we proceed analytically. We solve the SDP \eqref{eq: z2 SDP} with $\rho$ included analytically in the basis $\bb$ (cf.\  \S\ref{sec: parameter-dependent}) and with $U=\rho^4+u_0$, seeking a $V$ ansatz that yields $u_0^*=0$. When $r$ is fixed at 28 instead of being treated analytically, quartic $V$ apparently suffices to prove a sharp upper bound (cf.\ table~\ref{tab: enclosures}), so we first try a quartic ansatz for $V(x,y,z,\rho)$. The most general quartic $V$ we need to consider includes only symmetric terms, and its highest-degree terms must take the form $\rho^sx^p(y^2+z^2-2\rho z)^q$ in order to avoid any degree-5 terms in $\f\cdot\nabla V$. Numerical solution of the SDP gives $u_0^*\approx0$, suggesting that a general quartic $V$ is indeed sufficient. However, not all terms in this $V$ ansatz are needed. Through a process of trial and error we remove and combine terms while repeatedly computing $u^*$. If removing a term from $V$ makes $u^*$ exceed zero by more than numerical error, the term must be kept. The result of this process is a $V$ ansatz that has few tunable coefficients and thus leads to tractable analysis. In particular, the ansatz
\begin{equation}
V = c_1 \left[\tfrac{1}{\sigma}x^4+(y^2+z^2-2 \rho  z)^2 
			+ 8\rho ^2 \left(y^2+z^2-2 \rho  z\right) + \tfrac{6}{\sigma}\rho ^2 x^2\right] 
			 - c_2\rho  \left(\tfrac{1}{\sigma}x+y\right)^2
\label{eq: z3 V}
\end{equation}
gives $u^*\approx0$ at the standard values $(\beta,\sigma)=(8/3,10)$. The factors of $\frac{1}{\sigma}$ are included to avoid some parameter-dependence in $S_U$; terms in $V$ of the form $\tfrac{1}{\sigma}x^n$ produce $\sigma$-independent terms in $\f\cdot\nabla V$.

\subsubsection{\label{sec: z3 analytical}Analytical SDP solution}

The bound $\rho\ov{z^3}\le\rho^4$ will be proven if the coefficients $c_i$ in \eqref{eq: z3 V} can be chosen to make $S_U(x,y,z,\rho)$ an SOS polynomial when $U=\rho^4$. This SOS condition is enforced as in \eqref{eq: z2 SDP} by requiring that $S_U$ can be represented by Gram matrices $\Q_s,\Q_a\succeq0$. Applying the $V$ ansatz \eqref{eq: z3 V} to the definition \eqref{eq: S_U} of $S_U$ yields a homogenous quartic $S_U$ that can be represented using the quadratic monomial vectors
\begin{align}
\bb_s &= \begin{bmatrix} x^2 & xy & y^2 & \rho^2 & \rho z & z^2  \end{bmatrix}^T,
&
\bb_a &= \begin{bmatrix} \rho x & \rho y & xz & yz  \end{bmatrix}^T.
\label{eq: z3 basis monomial}
\end{align}
Any $S_U$ proving a sharp bound must vanish at the equilibria $\x^\pm$ that saturate the bound (cf.\ \S\ref{sec: why}), which is possible only if $S_U$ vanishes whenever $z=\rho$ or $x=y$. The above monomial basis vectors are thus unnecessarily general, and it suffices to represent $S_U$ using 
\begin{align}
\bhat_s &= \begin{bmatrix} x^2-xy \\ x^2-y^2 \\ \rho(z-\rho) \\ z(z-\rho) \end{bmatrix},
&
\bhat_a &= \begin{bmatrix} \rho(x-y) \\ x(z-\rho) \\ y(z-\rho) \end{bmatrix}.
\end{align}
Even with these smaller basis vectors, however, numerical SDP solutions give an approximate $\Qhat_s$ that is very close to singular, with a nearly zero eigenvalue corresponding to an eigenvector of approximately $\begin{bmatrix}0&0&1&1\end{bmatrix}^T$. 
This suggests that the exact $\Qhat_s$ has a null space spanned by this vector. We thus can reduce the dimension of $\bhat_s$ further by choosing
\begin{align}
\bhat_s &= \begin{bmatrix} x^2-xy \\ x^2-y^2 \\ (z-\rho)^2 \end{bmatrix},
&
\bhat_a &= \begin{bmatrix} \rho(x-y) \\ x(z-\rho) \\ y(z-\rho) \end{bmatrix},
\label{eq: z3 basis final}
\end{align}
which leads to a strictly feasible SDP.

The larger bases \eqref{eq: z3 basis monomial} were needed to solve the SDP optimization \eqref{eq: z2 SDP} that suggested the $V$ ansatz \eqref{eq: z3 V}, but the reduced bases \eqref{eq: z3 basis final} suffice to represent $S_U$ in the optimal case where $U=\rho^4$. While the reduction of basis vectors from \eqref{eq: z3 basis monomial} to \eqref{eq: z3 basis final} was guided by simple observations, more complicated cases may require computer-assisted reduction as in~\cite{Monniaux2011}. 

To complete the proof that $\rho\ov{z^3}\le\rho^4$ we must find symmetric $3\times3$ matrices $\Qhat_s,\Qhat_a\succeq0$ that represent $S_U$ with the basis vectors \eqref{eq: z3 basis final}. Matching coefficients of $S_U=\bhat_s^T\Qhat_s\bhat_s+\bhat_a^T\Qhat_a\bhat_a$ with those of the expression \eqref{eq: S_U} defining $S_U$ gives 12 independent affine relations between the coefficients of $V$ and the entries of $\Qhat_s$ and $\Qhat_a$. These relations uniquely determine the coefficients of $V$ as
\begin{align}
c_1&=\frac{1}{4\beta}, & c_2&=\frac{\sigma}{2(1+\sigma)},
\label{eq: z3 c_i}
\end{align}
and they determine the entries of the Gram matrices up to two degrees of freedom. Letting $\gamma_1$ and $\gamma_2$ denote the $(2,3)$ and $(3,3)$ entries of $\Qhat_a$, respectively, the Gram matrices must take the form
\begin{align}
\Qhat_s &= \tfrac{1}{2\beta}\begin{bmatrix}
2 & -1 & 1+2\beta\gamma_1 \\
-1 & 2 & \beta(\gamma_2-1)-1 \\
1+2\beta\gamma_1 & \beta(\gamma_2-1)-1 & 2\beta
\end{bmatrix},
\label{eq: z3 Qs}
\\
\Qhat_a &= \tfrac{1}{2}\begin{bmatrix}
\tfrac{6}{\beta} & -\tfrac{1}{1+\sigma} & -1 \\
-\tfrac{1}{1+\sigma} & 2(1-2\gamma_1-\gamma_2) & 2\gamma_1 \\
-1 & 2\gamma_1 & 2\gamma_2
\end{bmatrix}.
\label{eq: z3 Qa}
\end{align}

The matrices \eqref{eq: z3 Qs} and \eqref{eq: z3 Qa} represent $S_U$ for any $(\gamma_1,\gamma_2)$ but are not necessarily positive semidefinite. If it is possible to choose $(\gamma_1,\gamma_2)$ so that $\Qhat_s,\Qhat_a\succeq0$, then $S_U$ is SOS. Since the present SDP is fairly small, we proceed analytically to find the region of the $(\beta,\sigma)$ plane where the SDP is feasible -- that is, where there exist $(\gamma_1,\gamma_2)$ such that $\Qhat_s,\Qhat_a\succeq0$.

\begin{figure}[t]
\begin{center}
\includegraphics[width=240pt]{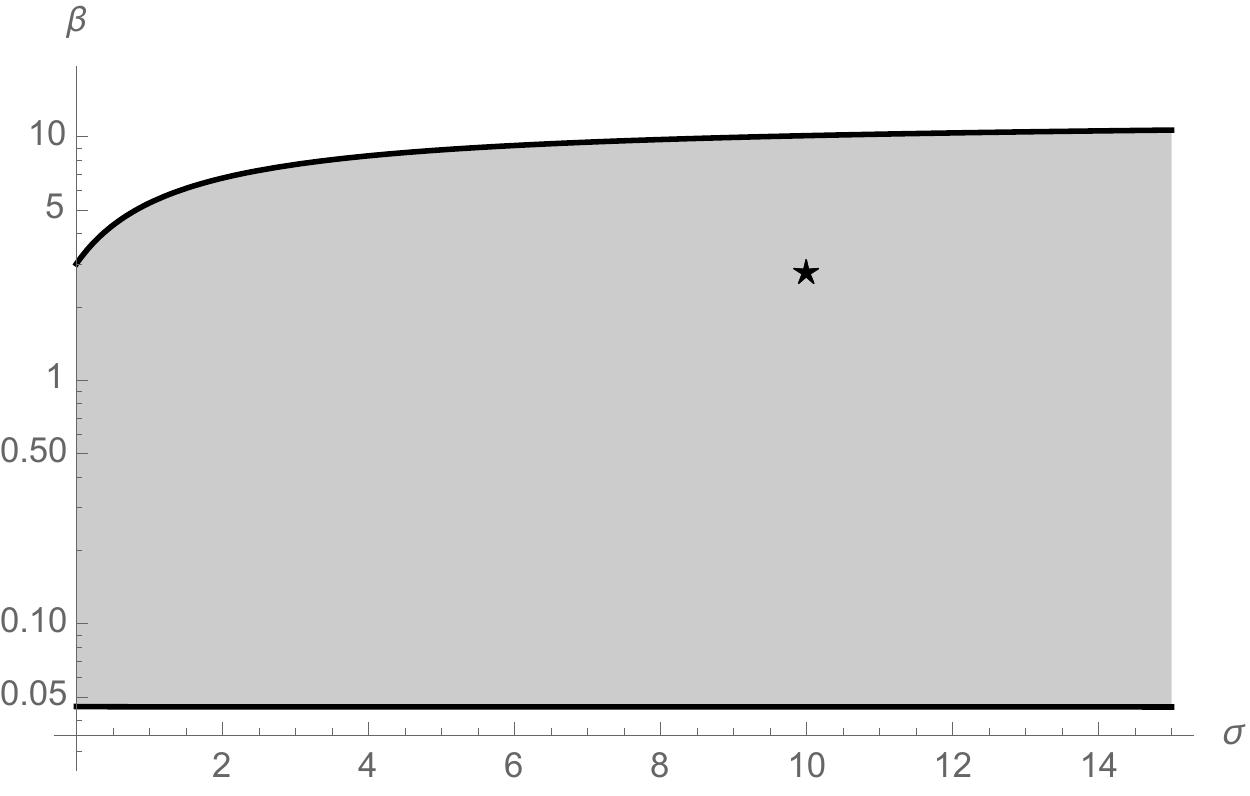}
\caption{Region of the $(\beta,\sigma)$ plane where we have proved that $\ov{z^3}\le(r-1)^3$ when $r\ge1$. The region includes the standard chaotic parameters ($\star$) and extends infinitely to the right as $\sigma\to\infty$. Outside the shaded region we have neither proved nor disproved the bound, but it cannot be proved using the $V$ ansatz \eqref{eq: z3 V}.}
\label{fig: z3}
\end{center}
\end{figure}

By Descartes' rule of signs \cite{Wang2004a}, a symmetric matrix is positive semidefinite if and only if its characteristic polynomial has coefficients that alternate between nonnegative and nonpositive. Requiring this of the characteristic polynomials of $\Qhat_s$ and $\Qhat_a$ gives six inequalities. One of these inequalities holds for all positive $\beta$, and the other five are
\beq
\begin{split}
- \beta^2\left[4 \gamma_1^2+2 (\gamma _2-1) \gamma _1 
 	+(\gamma _2-1)^2\right] + \beta(2+\gamma _2-2\gamma _1) - 1 &\ge0 \\
\beta^2\left[4\gamma _1^2+(\gamma _2-1)^2\right]
	+\beta\left[4 \gamma _1-2(\gamma _2+3)\right] - 1 &\le0 \\
2 \gamma _1 (\sigma +1) \left[\beta(\sigma +2)-12 \gamma _2 (\sigma +1)\right]
+\gamma _2 \left[(\beta +12) \sigma (\sigma +2)-12\gamma _2(\sigma+1)^2+12\right] \\
~-\beta(\sigma +1)^2 - 12\gamma _1^2 (\sigma +1)^2 &\ge0 \\
4 (\sigma +1)^2 \left[\beta\gamma _1^2+2 \gamma _1 (\beta  \gamma _2+3)
+\beta(\gamma _2-1) \gamma _2\right]+\beta \left[\sigma  (\sigma +2)+2\right]-12 (\sigma +1)^2 &\le0 \\
\beta(1-2\gamma _1)+3 &\ge0.
\end{split}
\label{eq: z3 ineqs}
\eeq
Whether there exist $(\gamma_1,\gamma_2)$ satisfying all five inequalities simultaneously depends on $(\beta,\sigma)$. The inequalities are too complicated to reduce by hand but are tractable using computer algebra. Quantifier elimination performed with the Mathematica syntax {\tt Reduce[Exists[\{$\gamma_1$,$\gamma_2$\},ineq]]}, where {\tt ineq} represents the inequalities \eqref{eq: z3 ineqs}, reveals that at each $\sigma>0$ there exist admissible $(\gamma_1,\gamma_2)$ over a bounded interval of positive $\beta$. Figure~\ref{fig: z3} shows part of the feasible region in the $(\beta,\sigma)$ plane. The upper extent of $\beta$ is
\[ \beta \le 12\frac{1+2\sigma+\sigma^2}{4+4\sigma+\sigma^2}, \]
which increases from 3 to 12 as $\sigma$ varies from 0 to $\infty$. The lower extent of $\beta$ is the smallest root of a particular degree-10 polynomial (omitted for brevity) whose coefficients are polynomials in $\sigma$; this root decreases very slightly from approximately $0.0456122$ to $0.0454294$ as $\sigma$ varies from 0 to $\infty$.

Since the shaded region of the $(\beta,\sigma)$ plane in figure~\ref{fig: z3} is where $(\gamma_1, \gamma_2)$ can be chosen so that $\Qhat_s,\Qhat_a\succeq0$, these are the parameter values where we have proved that $\rho\ov{z^3}\le\rho^4$. That is, $\ov{z^3}\le(r-1)^3$ on every trajectory of the Lorenz system for these $(\beta,\sigma)$ and all $r\ge1$. We cannot say whether the same result holds outside the shaded region. If so, a $V$ ansatz other than \eqref{eq: z3 V} is needed to prove the bound because for this ansatz the polynomial $S_U$ cannot be SOS.

\subsubsection{\label{sec: z3 explicit}Explicit certificate at the standard parameters}

The preceding proof that $\rho\ov{z^3}\le\rho^4$ in the $(\beta,\sigma)$ regime of figure~\ref{fig: z3} does not provide an easily checked ``certificate'' that $S_U$ is indeed SOS. Computer algebra was used to determine the $(\beta,\sigma)$ at which there exist $(\gamma_1,\gamma_2)$ making the Gram matrices \eqref{eq: z3 Qs} and \eqref{eq: z3 Qa} positive semidefinite, but no explicit expressions were found for these (non-unique) Gram matrices. Here we illustrate how to find explicit Gram matrices at the standard values of $(\beta,\sigma)$ by selecting an admissible pair $(\gamma_1,\gamma_2)$.

At the standard values $(\beta,\sigma)=(3/8,10)$, the inequalities \eqref{eq: z3 ineqs} that are equivalent to the conditions $\Qhat_s,\Qhat_a\succeq0$ reduce to the four inequalities
\begin{align}
\frac{1}{33} (8-33 \gamma_2)-\frac{1}{33} \sqrt{89} \sqrt{9 \gamma_2-2} &\leq
	\gamma_1\leq\frac{1}{33} (8-33 \gamma_2)
\label{eq: standard gamma cond 1}
		+\frac{1}{33} \sqrt{89} \sqrt{9 \gamma_2-2} \\
\frac{2}{9} &\le \gamma_2 \le 0.2370\ldots,
\label{eq: standard gamma cond 2}
\end{align}
where $0.2370\ldots$ is the smallest real root of the quartic polynomial
\[ 909988849 - 9432927504 \lambda + 49128348096 \lambda^2 - 118056102912 \lambda^3 + 43717791744 \lambda^4. \]
The admissible set of $(\gamma_1,\gamma_2)$ is highly constrained but still two-dimensional. One admissible point that is particularly simple is $(\gamma_1,\gamma_2)=(0,3/8)$. Applying this choice to the Gram matrices \eqref{eq: z3 Qs} and \eqref{eq: z3 Qa} gives an explicit representation of $S_U$ by positive definite matrices,
\begin{equation}
S_U(x,y,z,\rho) = \bhat_s^T\begin{bmatrix}
\tfrac{3}{8} & -\tfrac{3}{16} & \tfrac{3}{16} \\[3pt]
-\tfrac{3}{16} & \tfrac{3}{8} & -\frac{1}{2} \\[3pt]
\tfrac{3}{16} & -\tfrac{1}{2} & 1 \end{bmatrix} \bhat_s +
\bhat_a^T\begin{bmatrix}
\tfrac{9}{8} & -\tfrac{1}{22} & -\tfrac{1}{2} \\[3pt]
-\tfrac{1}{22} & \tfrac{5}{8} & 0 \\[3pt]
-\tfrac{1}{2} & 0 & \tfrac{3}{8}
\end{bmatrix} \bhat_a,
\label{eq: z3 Q standard}
\end{equation}
where the basis vectors $\bhat_s$ and $\bhat_a$ are as defined by \eqref{eq: z3 basis final}.

The entire proof that $\rho\ov{z^3}\le\rho^4$ when $(\beta,\sigma)=(8/3,10)$ thus can be summarized as follows. For any differentiable $V$,
\beq
\rho\ov{z^3} = \rho^4 - \ov{S_U},
\label{eq: z3 explicit}
\eeq
where in this case $S_U=-\left[\rho z^3 - \rho^4 + \f\cdot\nabla V\right]$. Let $V$ be defined by the ansatz \eqref{eq: z3 V} with the coefficients \eqref{eq: z3 c_i}. For this $V$, one can check by exact arithmetic that $S_U$ is indeed equal to expression \eqref{eq: z3 Q standard}, and that the matrices in that expression are positive definite. Thus $\ov{S_U}$ in \eqref{eq: z3 explicit} is nonnegative, and the bound $\rho\ov{z^3}\le\rho^4$ is proven.

Verifying that the matrices in \eqref{eq: z3 Q standard} are positive definite confirms that $S_U$ is SOS. Nothing further is proved by finding an explicit SOS representation, but one can be found if desired by computing the Cholesky decompositions $\Qhat_s=L_sL_s^T$ and $\Qhat_a=L_aL_a^T$. This gives the SOS expression $S_U = \Vert L_s^T\bhat_s\Vert_2^2 + \Vert L_a^T\bhat_a\Vert_2^2$. For the Gram matrices in \eqref{eq: z3 Q standard} the Cholesky decompositions are
\begin{align}
L_s^T &= \tfrac{1}{\sqrt2}\begin{bmatrix}
\tfrac{\sqrt{3}}{2} & -\tfrac{\sqrt{3}}{4} & \tfrac{\sqrt{3}}{4} \\[4pt]
0 & \tfrac{3}{4} & -\tfrac{13}{12} \\[4pt]
0 & 0 & \tfrac{\sqrt{23}}{6} \end{bmatrix},
&
L_a^T &= \tfrac{1}{\sqrt2}\begin{bmatrix}
\tfrac{3}{2} & -\tfrac{2}{33} & -\tfrac{2}{3} \\[4pt]
0 & \tfrac{\sqrt{5429}}{66} & -\tfrac{8}{3\sqrt{5429}} \\[4pt]
0 & 0 & \tfrac{\sqrt{6607}}{2\sqrt{5429}}
\end{bmatrix}.
\end{align}
The equality $S_U = \Vert L_s^T\bhat_s\Vert_2^2 +  \Vert L_a^T\bhat_a\Vert_2^2$, after some slight simplification in the first term, gives the SOS representation
\begin{multline}
S_U(x,y,z,\rho) = \tfrac{3}{32}\left[ (x-y)^2 + (z-\rho)^2 \right]^2 
	+ \tfrac{1}{288}\left[ 9(x^2 - y^2) - 13(z-\rho)^2 \right]^2 \\
	+ \tfrac{23}{72}(z-\rho)^2
	+ \tfrac{1}{2}\left[ \tfrac{3}{2}\rho(x-y) - \tfrac{2}{33}x(z-\rho) - \tfrac{2}{3}y(z-\rho) \right]^2 \\
	+ \tfrac{1}{2}\left[ \tfrac{\sqrt{5429}}{66}x(z-\rho) -\tfrac{8}{3\sqrt{5429}}y(z-\rho) \right]^2
	+ \tfrac{6607}{21716}y^2(z-\rho)^2.
\end{multline}

\subsection{\label{sec: lower}Lower bounds of zero in the Lorenz system}

Averages of all symmetric moments of the Lorenz system up to quartic degree are nonnegative at the standard parameters, meaning they are minimized by trajectories on or approaching the zero equilibrium. Since symmetric moments are those taking the form $\mom$ with $l+m$ even, the exponents $l$ and $m$ are either both even or both odd. When both are even, $\mom\ge0$ holds not only on average but also everywhere on the global attractor since $z\ge0$ on the attractor at all positive parameters. On the other hand, the five moments with $l$ and $m$ odd ($xy$, $xyz$, $x^3y$, $xy^3$, $xyz^2$) are negative on parts of the attractor, at least at the standard parameters. Nonetheless, the time averages of these moments are nonnegative along every trajectory. For the four moments other than $xy^3$ these lower bounds follow easily at all positive parameters from \eqref{eq: proportional}: the average of each moment is proportional to the average of a different moment that is nonnegative everywhere on the attractor. Below we prove that $\ov{xy^3}\ge0$ when $r\ge0$ and $\beta\in\left[6-4\sqrt{2},6+4\sqrt{2}\right]\approx[0.34,11.66]$, which includes the standard value $\beta=8/3$. It is an open question whether averages of higher-degree moments with odd $l$ and $m$ are all nonnegative also.

To determine a $V$ ansatz that suffices to prove $\ov{xy^3}\ge0$, we numerically solve the SDP optimization \eqref{eq: L optimal Lorenz} with $r$ included analytically in the basis $\bb$. The analytical treatment of $r$ precludes showing $\ov{xy^3}\ge0$ directly since the result is false for negative $r$, but we can instead show that $r\ov{xy^3}\ge0$ for all $r$. We thus let $\varphi=rxy^3$ and search for a $V$ that give $L^*=0$. At the standard values $(\beta,\sigma)=(8/3,10)$ we find $L^*\approx0$ using a general quartic ansatz for $V$. Further numerical trial-and-error suggests that only four terms from the $V$ basis are needed, and the coefficients happen to be independent of $\beta$ and~$\sigma$:
\beq
V(y,z,r) = -r^2z^2 + ry^2z + \tfrac{4}{3}rz^3 - \tfrac{1}{2}(y^2+z^2)^2.
\label{eq: xy3 V}
\eeq
Applying the above $V$ to the definition \eqref{eq: S_L} of $S_L$ gives a polynomial independent of $x$,
\beq
S_L(y,z,r) = 2\beta r^2 z^2-(2+\beta) r y^2 z-
			4 \beta r z^3 + 2y^4 + 2(1+\beta)y^2z^2 + 2\beta z^4.
\label{eq: xy3 S}
\eeq
The bound $r\ov{xy^3}\ge0$ will be proven if we can show that $S_L$ is SOS, or equivalently that it can be represented by some $\Q_s,\Q_a\succeq0$.

The homogenous quartic polynomial $S_L$ can be represented using the basis vectors
\begin{align}
\bb_s&=\begin{bmatrix}rz&y^2&z^2\end{bmatrix}^T,
&
\bb_a&= \begin{bmatrix}yz\end{bmatrix}.
\end{align}
We must choose a $3\times3$ symmetric $\Q_s$ and a scalar $\Q_a$ such that $S_L=\bb_s^T\Q_s\bb_s+\Q_ay^2z^2$. Matching the coefficients of this expression with those of \eqref{eq: xy3 S} gives six constraints that let all entries of $\Q_s$ be either fixed or expressed in terms of $\Q_a$,
\begin{equation}
Q_s=\begin{bmatrix}
4\beta & -(2+\beta) & -4\beta \\
-(2+\beta) & 4 & \tfrac{1}{2}(4+4\beta-\Q_a) \\
-4\beta & \tfrac{1}{2}(4+4\beta-\Q_a) & 4\beta
\end{bmatrix}.
\label{eq: xy3 Qs}
\end{equation}
The polynomial $S_L$ is SOS if and only if we can choose $\Q_a\ge0$ such that $\Q_s\succeq0$. Whether this is possible depends on the value of $\beta$. The characteristic polynomial of $Q_s$ takes the form $\lambda^3 + c_2\lambda^2 + c_1\lambda+c_0=0$, and by Descartes' rule of signs its roots are all nonnegative if and only if $c_2\le0$, $c_1\ge0$, and $c_0\le0$. The inequality $c_2\le0$ holds for all positive $\beta$ since $c_2=-4(1+2\beta)$. The inequality $c_0\le0$ requires that $\Q_a=2\beta$ since $c_0=\beta(\Q_a-2\beta)^2$. For this $\Q_a$ the remaining condition is $c_1=-(4-12\beta+\beta^2)\ge0$, which holds if and only if $\beta\in\left[6-4\sqrt{2},6+4\sqrt{2}\right]$. For $\beta$ in this range the polynomial $S_L$ given by \eqref{eq: xy3 S} is SOS, thereby proving that $r\ov{xy^3}\ge0$. Numerical SDP solutions suggest that the bound holds for other $\beta$ values also, but proving this would require a different choice of $V$.

Finding an explicit SOS representation of $S_L$ does not prove anything further but can be done easily if desired. The matrix \eqref{eq: xy3 Qs} with $\Q_a=2\beta$ factors into $Q_s=L_sL_s^T$, where
\beq
L_s^T=\begin{bmatrix}
2\sqrt{\beta} &  -\frac{2+\beta}{2\sqrt{\beta}} & -2\sqrt{\beta} \\
0 & \frac{\sqrt{-(\beta^2-12\beta+4)}}{2\sqrt\beta} & 0
\end{bmatrix}.
\eeq
The polynomial $S_L$ given by \eqref{eq: xy3 S} therefore can be written as
\begin{align}
S_L(r,y,z)	&= \Vert L_s^T\bb_s\Vert_2^2 + 2\beta y^2z^2 \\
	&= \tfrac{1}{4\beta}\left[4\beta rz -(2+\beta)y^2-4\beta z^2\right]^2
		+\tfrac{-(\beta ^2-12 \beta+4)}{4\beta}y^4+2 \beta  y^2 z^2.
\label{eq: xy3 S sos}
\end{align}
Expression \eqref{eq: xy3 S sos} is SOS when the coefficient of $y^4$ is nonnegative. This occurs if and only if $\beta\in\left[6-4\sqrt{2},6+4\sqrt{2}\right]$, which is exactly the condition for $\Q_s\succeq0$ and for $L_s$ to be real.

\section{\label{sec: con}Conclusions}

We have described methods for bounding time averages in dynamical systems using semidefinite programming, and we have applied these methods to the Lorenz system. The bounds are global in the sense that they apply to all possible trajectories. Rigorous bounds have been obtained in two different ways: by analytically finding exact solutions to SDPs (or showing that such solutions exist), and by using interval arithmetic to enclose approximate numerical solutions. The former method is needed to prove sharp bounds, while the latter method is easier to implement and can produce nearly sharp bounds. Most of the bounds we have constructed are novel, many are very tight, and some are perfectly sharp, thereby demonstrating that the complicated phase space of a chaotic system does not prevent the SDP approach from succeeding. We are not aware of any other way to produce rigorous bounds of this quality, except in simple cases where sharp bounds can be constructed without computer assistance.

\begin{table}[t]
\begin{center}
\caption{\label{tab: summary}Time averages and our best upper bounds for moments of the Lorenz system, normalized by each moment's value \eqref{eq: mom nonzero} at the nonzero equilibria. Also shown is the tightness of each bound as measured by its margin above the maximum known average. For each moment the maximum known average occurs either on the nonzero equilibria (where its normalized value is 1) or on the $\oplus$$\ominus$ periodic orbit. Moments grouped together (e.g.\ $z$, $x^2$, $xy$) have identical normalized averages according to \eqref{eq: proportional}. Chaotic averages are obtained by numerical integration (cf.\ \S\ref{sec: extremal traj}), and most are converged to the precision shown, but for the $z$-independent quartic moments the last digit is uncertain.}
\begin{tabular}{l l l l l}
Moment & Chaotic mean & Maximum mean
	& Upper bound & Tightness \\
	\hline
$z$, $x^2$, $xy$	& 0.87223	& 1 				& 1	& sharp \\ 
$y^2$ 			& 1.12780	& 1.1621684		& 1.1627	& 0.046\% \\ 
$z^2$, $xyz$ 		& 0.86276 	& 1 			& 1	& sharp \\ 
$x^2z$ 			& 0.96689	& 1 				& 1.0000003	& 0.00003\% \\ 
$y^2z$ 			& 1.02733	& 1.0394975 		& 1.0396	& 0.0099\% \\ 
$z^3$, $xyz^2$ 	& 0.93716 	& 1 			& 1	& sharp \\	
$x^4$, $x^3y$		& 1.74779	& 1.9111906 		& 1.9164	& 0.27\% \\ 
$x^2y^2$			& 2.07089	& 2.2975630 		& 2.3220	& 1.06\% \\ 
$x^2z^2$			& 1.15101	& 1.1893425 		& 1.1899	& 0.047\% \\ 
$xy^3$  			& 2.65721	& 2.9987454 		& 3.0239	& 0.84\% \\ 
$y^4$ 			& 3.62466 	& 4.1459937 	& 4.1842	& 0.92\% \\ 
$y^2z^2$			& 1.03615	& 1.0484088 		& 1.0489	& 0.047\% \\  
$z^4$			& 1.09006	& 1.1155092		& 1.1158	& 0.026\%	 
\end{tabular}
\end{center}
\end{table}

The quantities we have bounded in the Lorenz system are the eighteen time-averaged moments $\ov\mom$ up to quartic degree that are invariant under the symmetry of the system. Sharp bounds are possible when the trajectories that saturate them lie on (or tend to) equilibrium points. Sharp lower bounds of zero hold for all eighteen moments, and they are saturated by the zero equilibrium. The only one of these bounds that is nontrivial to prove is $\ov{xy^3}\ge0$. Our best upper bounds at the standard parameter values are collected in table~\ref{tab: summary}, along with each moment's chaotic average and maximum known average. Sharp upper bounds have been proved for seven moments, and these are saturated by the two nonzero equilibria. All seven of these results can be derived from the bound $\ov{z^3}\le(r-1)^3$ that we have proved for $r\ge1$ and a range of $(\beta,\sigma)$. The upper bounds on the remaining eleven moments, computed using interval arithmetic with the software VSDP, are all within 1\% of being sharp.

One remaining challenge is to construct sharp bounds when this requires exact optimal solutions to SDPs that are too large to study analytically. Such an example is the conjectured bound $\ov{x^2z}\le\beta(r-1)^2$ in the Lorenz system at the standard parameters; our upper bound on $\ov{x^2z}$ was constructed using interval arithmetic and so is very slightly conservative. The sharp bound might be proved using symbolic-numerical algorithms that project approximate numerical solutions onto exact rational ones \cite{Peyrl2008, Kaltofen2012, Wu2013}. A potential difficulty is that these algorithms are guaranteed to succeed (given enough precision) only if the SDP to be solved is strictly feasible, whereas exactly optimal solutions are marginally feasible in general. Projection might still succeed in marginal cases \cite{Peyrl2008}, and if not one can formulate a strictly feasible SDP by reducing the polynomial basis. We have carried out such reduction analytically in \S\ref{sec: z and z2} and \S\ref{sec: z3}, and the algorithm of \cite{Monniaux2011} offers an automated approach for larger SDPs. There are cases where a Gram matrix with rational entries may not exist, however, including when the desired bound is irrational.

Bounds that apply to all trajectories of a dynamical system yield global information in a way that computing particular trajectories cannot. Computing every periodic orbit, for instance, is impossible when there are infinitely many. On the other hand, global bounds can be unnecessarily conservative if one is interested only in particular trajectories. In the Lorenz system, one might want bounds that apply to chaotic trajectories but not necessarily to unstable periodic orbits or equilibria. The nonzero equilibria are separate from the strange attractor, and bounds that do not apply to these equilibria could be proved by enforcing bounding conditions only on a subset of phase space that omits them. This approach succeeded for the van der Pol oscillator \cite{Fantuzzi2016}, but our preliminary efforts with the Lorenz system have been plagued by poor numerical conditioning. We are not aware of a method for excluding unstable trajectories that are embedded in the strange attractor, such as periodic orbits or the equilibrium at the origin. Adding finite noise as in \cite{Fantuzzi2016} can give bounds on stochastic expectations that are close to chaotic time averages, but here too we have had numerical difficulties. Proving bounds that apply only to particular trajectories seems to require progress both in formulating bounding conditions as SDPs and in solving ill-conditioned SDPs. Nevertheless, the methods used in the present work can produce very tight global bounds in chaotic systems. Applying them to nearly any low-dimensional dissipative system is likely to yield novel results.

\subsection*{Acknowledgements}
The author was partially supported during this work by the James Van Loo Post-Doctoral Fellowship and NSF award DMS--1515161. The author thanks Divakar Viswanath for providing numerically computed periodic orbits of the Lorenz system. This work has benefited from conversations with Sergei Chernyshenko, Charles R.\ Doering, Bruno Eckhardt, Giovanni Fantuzzi, and Anthony Quas.

\appendix

\section{\label{app: vsdp}Details of verified computations using VSDP}

When computing the bounds reported in tables \ref{tab: enclosures} and \ref{tab: non-sharp} with the software VSDP, we most often rescaled the Lorenz system by $\x\mapsto20\x$ and then converted bounds back to the original scaling. Various researchers using SOS methods to study dynamics (although not to bound time averages) have found that rescaling the governing equations improves the convergence of SDP solvers. A heuristic that is often used (e.g.\ \cite{Henrion2014}) is to rescale each coordinate of $\x\in\mathbb R^n$ so that the relevant dynamics occur approximately in the cube $[-1,1]^n$. Our rescaling is similar, putting most of the Lorenz attractor in the domain $[-1,1]^2\times[0,2]$. We find that rescaling by 10 or 40 instead of 20 usually gives more conservative upper enclosures $U^+$. Results become more sensitive to the rescaling when the auxiliary function $V$ reaches degree 8 or 10, at which point slightly different scalings can significantly affect the conservativeness of the enclosures. The degree-10 bounds in table~\ref{tab: non-sharp} were produced after rescaling by either 25 (for $x^4$, $x^2y^2$, $xy^3$) or 30, and the degree-8 lower enclosure in table~\ref{tab: enclosures} was produced after rescaling by 10. The time required to compute each tabulated bound on a single processor ranged from seconds to minutes.

Results reported here for $V$ of degree 6, 8, or 10 are not quite rigorous to the standard of a computer-assisted proof. This is because we have used the software YALMIP \cite{Lofberg2004, Lofberg2009} to automatically parse the SOS conditions, formulate corresponding SDPs, and pass them to the VSDP software. This incurs roundoff errors that are not accounted for since YALMIP does not use interval arithmetic. A parser that uses interval arithmetic is under development. Until it is available, rigorous results can be found by formulating the relevant SDPs manually. We have done this only for $V$ of degree 4, in which cases roundoff errors introduced by YALMIP have all been orders of magnitude smaller than the margins of the enclosures generated by VSDP.

\vspace{-4pt}

\section{\label{app: equalities}Relations between mean moments}

\begin{table}[t]
\caption{\label{tab: all relations}Symmetric mean moments up to quartic degree expressed as linear combinations of moments that either are in the minimal set $\big\{ \ov z, \ov{z^2}, \ov{y^2z}, \ov{z^3}, \ov{y^2z^2}, \ov{z^4} \big\}$ or are expressed in this same way higher in the table. Each equality is an instance of the identity $\ov\varphi=\ov{\varphi+\f\cdot\nabla V}$ with $\varphi$ and $V$ defined as shown.}
\begin{align*}
\ov\varphi	&= \ov{\varphi+\f\cdot\nabla V}			 	&& V(x,y,z) \\[2pt]
\hline
\ov{xy}		&=\beta\ov{z}		&& -z \\
\ov{x^2}		&=\ov{xy} 			&& \tfrac{1}{2\sigma}x^2 \\
\ov{y^2}		&=r \ov{xy}-\beta\ov{z^2} 	&& \tfrac12(y^2 + z^2) \\
\ov{x^2z}		&=r\ov{x^2} - (1+\sigma)\ov{xy} + \sigma\ov{y^2}	&& xy \\
\ov{xyz}		&=\beta \ov{z^2}	&& -\tfrac12 z^2 \\
\ov{x^3y}		&=(\beta+2\sigma)\ov{x^2z} - 2\sigma\ov{xyz} && -x^2z \\
\ov{x^4} 		&= \ov{x^3y}	&& \tfrac{1}{4\sigma}x^4 \\
\ov{xyz^2} 	&= \beta\ov{z^3}	&& -\tfrac13z^3 \\
\ov{xy^3}		&= (2+\beta)\ov{y^2z}-2r\ov{xyz}+2\ov{xyz^2} && -y^2z \\
\ov{x^2z^2}	&= \frac{1}{1+\beta+2\sigma}\Big[ 
					\sigma(\sigma+1)\ov{y^2z} + r(1+\sigma)\ov{x^2z}
					&& \frac{2(1+\sigma)xyz + x^2(y^2+z^2)}
						{2(1+\beta+2\sigma)} \\
				& \hspace{18pt} - (1+\sigma)(1+\beta+\sigma)\ov{xyz}  + r\ov{x^3y} 
					+ \sigma\ov{xy^3} + \sigma\ov{xyz^2} \Big] \\
\ov{x^2y^2}	&= -\sigma\ov{y^2z}-r\ov{x^2z}+(1+\beta+\sigma)\ov{xyz}
					+\ov{x^2z^2} && -xyz \\
\ov{y^4}		&= r\ov{xy^3} + r\ov{xyz^2} - (1+\beta)\ov{y^2z^2} - \beta\ov{z^4}
					&& \tfrac14(y^2+z^2)^2 \\
\end{align*}
\end{table}

On any trajectory of the Lorenz system, all mean moments up to degree 2, 3, or 4 that are symmetric under $(x,y)\mapsto(-x,-y)$ can be expressed as linear combinations of $\big\{ \ov z, \ov{z^2} \big\}$, $\big\{ \ov z, \ov{z^2}, \ov{y^2z}, \ov{z^3} \big\}$, or $\big\{ \ov z, \ov{z^2}, \ov{y^2z}, \ov{z^3}, \ov{y^2z^2}, \ov{z^4} \big\}$, respectively. Many of these relations are not useful for constructing bounds, although they are derived using the same basic identity $\ov{\f\cdot\nabla V}=0$ that is central to our bounding methods. When proving bounds we have sought $V$ such that $\varphi+\f\cdot\nabla V$ obeys the desired bound at all points in phase space. For the different objective of expressing various moments in terms of a smaller subset of moments, $V$ must be chosen so that $\varphi+\f\cdot\nabla V$ contains only moments from this subset. Table~\ref{tab: all relations} gives relations for symmetric moments of the Lorenz system, as well as the choices of $V$ that yield these relations. About half of these relations have appeared in the literature for decades \cite{Lucke1976a}. Every tabulated moment is expressed as a linear combination of moments that either are in the minimal set $\big\{ \ov z, \ov{z^2}, \ov{y^2z}, \ov{z^3}, \ov{y^2z^2}, \ov{z^4} \big\}$ or are expressed in this same way higher in the table. Each moment can be expressed using only the minimal set after some further substitution: the $\ov{x^2}$ relation and the one above it give $\ov{x^2}=\beta\ov z$, the $\ov{y^2}$ relation and the ones above it give $\ov{y^2}=\beta\big(r\ov z-\ov{z^2}\big)$, the $\ov{x^2z}$ relation and the ones above it give $\ov{x^2z}=\beta(1+\sigma)(r-1)\ov z-\beta\sigma\ov{z^2}$, and so on.

\bibliography{manuscript.bbl}

\end{document}